\documentclass{article}

\pdfoutput=1

\usepackage{etex}
\usepackage{tikz-cd}

\newenvironment{tikzar}[1]{\begin{tikzcd}[ampersand replacement=\&,#1]}%
{\end{tikzcd}}
\tikzset{
  no line/.style={draw=none,
    commutative diagrams/every label/.append style={/tikz/auto=false}
  },
}
\makeatletter

\def\tikzcdset{\pgfqkeys{/tikz/commutative diagrams}}
\def\tikzcd@setanchor#1[#2]#3\relax{%
  \ifx\relax#2\relax\else%
    \tikzcdset{@#1transform/.append style={#2},@shiftabletopath}%
  \fi%
  \ifx\relax#3\relax%
    \pgfutil@namelet{tikzcd@#1anchor}{pgfutil@empty}%
  \else%
    \pgfutil@namedef{tikzcd@#1anchor}{.#3}%
  \fi}

\tikzcdset{
  @shiftabletopath/.style={
    /tikz/execute at begin to={%
      \begingroup%
        \def\tikz@tonodes{coordinate[pos=0,commutative diagrams/@starttransform/.try](tikzcd@nodea) %
                          coordinate[pos=1,commutative diagrams/@endtransform/.try](tikzcd@nodeb)}%
        \path (\tikztostart) \tikz@to@path;%
      \endgroup%
      \def\tikztostart{tikzcd@nodea}%
      \tikzset{insert path={(tikzcd@nodea)}}},
    /tikz/commutative diagrams/@shiftabletopath/.code={}},
  start anchor/.code={%
    \pgfutil@ifnextchar[{\tikzcd@setanchor{start}}{\tikzcd@setanchor{start}[]}#1\relax},
  end anchor/.code={%
    \pgfutil@ifnextchar[{\tikzcd@setanchor{end}}{\tikzcd@setanchor{end}[]}#1\relax},
  start anchor=,
  end anchor=,
  shift left/.style={
    /tikz/commutative diagrams/@shiftabletopath,
    /tikz/execute at begin to={%
      \pgfpointnormalised{%
        \pgfpointdiff{\pgfpointanchor{tikzcd@nodeb}{center}}{\pgfpointanchor{tikzcd@nodea}{center}}}%
      \pgfgetlastxy{\tikzcd@x}{\tikzcd@y}%
      \pgfmathparse{#1}%
      \coordinate (tikzcd@nodea) at ([shift={(\pgfmathresult*\tikzcd@y,-\pgfmathresult*\tikzcd@x)}]tikzcd@nodea);%
      \coordinate (tikzcd@nodeb) at ([shift={(\pgfmathresult*\tikzcd@y,-\pgfmathresult*\tikzcd@x)}]tikzcd@nodeb);%
      \tikzset{insert path={(tikzcd@nodea)}}}},
  shift right/.style={
    /tikz/commutative diagrams/shift left={-(#1)}},
  transform nodes/.style={
    /tikz/commutative diagrams/@shiftabletopath,
    /tikz/commutative diagrams/@starttransform/.append style={#1},
    /tikz/commutative diagrams/@endtransform/.append style={#1}},
  shift/.style={
    /tikz/shift={#1},
    /tikz/commutative diagrams/transform nodes={/tikz/shift={#1}}},
  xshift/.style={
    /tikz/xshift={#1},
    /tikz/commutative diagrams/transform nodes={/tikz/xshift={#1}}},
  yshift/.style={
    /tikz/yshift={#1},
    /tikz/commutative diagrams/transform nodes={/tikz/yshift={#1}}}}

\makeatother

\usepackage{amssymb}
\usepackage{amsmath}
\usepackage{amsthm}
\usepackage{txfonts}
\usepackage{bbm}
\usepackage{enumerate}
\usepackage{fullpage}

\usepackage[T1]{fontenc}
\usepackage[utf8]{inputenc}
\usepackage{authblk}

\theoremstyle{plain}
\newtheorem{theorem}{Theorem}
\newtheorem{proposition}[theorem]{Proposition}
\newtheorem{lemma}[theorem]{Lemma}
\newtheorem{corollary}[theorem]{Corollary}
\newtheorem{definition}[theorem]{Definition}

\newdimen\proofrulebreadth \proofrulebreadth=.05em
\newdimen\proofdotseparation \proofdotseparation=1.25ex
\newdimen\proofrulebaseline \proofrulebaseline=2ex
\newcount\proofdotnumber \proofdotnumber=3
\let\then\relax
\def\hfi{\hskip0pt plus.0001fil}
\mathchardef\squigto="3A3B
\newif\ifinsideprooftree\insideprooftreefalse
\newif\ifonleftofproofrule\onleftofproofrulefalse
\newif\ifproofdots\proofdotsfalse
\newif\ifdoubleproof\doubleprooffalse
\let\wereinproofbit\relax
\newdimen\shortenproofleft
\newdimen\shortenproofright
\newdimen\proofbelowshift
\newbox\proofabove
\newbox\proofbelow
\newbox\proofrulename
\def\shiftproofbelow{\let\next\relax\afterassignment\setshiftproofbelow\dimen0 }
\def\shiftproofbelowneg{\def\next{\multiply\dimen0 by-1 }%
\afterassignment\setshiftproofbelow\dimen0 }
\def\setshiftproofbelow{\next\proofbelowshift=\dimen0 }
\def\setproofrulebreadth{\proofrulebreadth}

\def\prooftree{
\ifnum  \lastpenalty=1
\then   \unpenalty
\else   \onleftofproofrulefalse
\fi
\ifonleftofproofrule
\else   \ifinsideprooftree
        \then   \hskip.5em plus1fil
        \fi
\fi
\bgroup
\setbox\proofbelow=\hbox{}\setbox\proofrulename=\hbox{}%
\let\justifies\proofover\let\leadsto\proofoverdots\let\Justifies\proofoverdbl
\let\using\proofusing\let\[\prooftree
\ifinsideprooftree\let\]\endprooftree\fi
\proofdotsfalse\doubleprooffalse
\let\thickness\setproofrulebreadth
\let\shiftright\shiftproofbelow \let\shift\shiftproofbelow
\let\shiftleft\shiftproofbelowneg
\let\ifwasinsideprooftree\ifinsideprooftree
\insideprooftreetrue
\setbox\proofabove=\hbox\bgroup$\displaystyle 
\let\wereinproofbit\prooftree
\shortenproofleft=0pt \shortenproofright=0pt \proofbelowshift=0pt
\onleftofproofruletrue\penalty1
}

\def\eproofbit{
\ifx    \wereinproofbit\prooftree
\then   \ifcase \lastpenalty
        \then   \shortenproofright=0pt  
        \or     \unpenalty\hfil         
        \or     \unpenalty\unskip       
        \else   \shortenproofright=0pt  
        \fi
\fi
\global\dimen0=\shortenproofleft
\global\dimen1=\shortenproofright
\global\dimen2=\proofrulebreadth
\global\dimen3=\proofbelowshift
\global\dimen4=\proofdotseparation
\global\count255=\proofdotnumber
$\egroup  
\shortenproofleft=\dimen0
\shortenproofright=\dimen1
\proofrulebreadth=\dimen2
\proofbelowshift=\dimen3
\proofdotseparation=\dimen4
\proofdotnumber=\count255
}

\def\proofover{
\eproofbit 
\setbox\proofbelow=\hbox\bgroup 
\let\wereinproofbit\proofover
$\displaystyle
}%
\def\proofoverdbl{
\eproofbit 
\doubleprooftrue
\setbox\proofbelow=\hbox\bgroup 
\let\wereinproofbit\proofoverdbl
$\displaystyle
}%
\def\proofoverdots{
\eproofbit 
\proofdotstrue
\setbox\proofbelow=\hbox\bgroup 
\let\wereinproofbit\proofoverdots
$\displaystyle
}%
\def\proofusing{
\eproofbit 
\setbox\proofrulename=\hbox\bgroup 
\let\wereinproofbit\proofusing
\kern0.3em$
}

\def\endprooftree{
\eproofbit 
  \dimen5 =0pt
\dimen0=\wd\proofabove \advance\dimen0-\shortenproofleft
\advance\dimen0-\shortenproofright
\dimen1=.5\dimen0 \advance\dimen1-.5\wd\proofbelow
\dimen4=\dimen1
\advance\dimen1\proofbelowshift \advance\dimen4-\proofbelowshift
\ifdim  \dimen1<0pt
\then   \advance\shortenproofleft\dimen1
        \advance\dimen0-\dimen1
        \dimen1=0pt
        \ifdim  \shortenproofleft<0pt
        \then   \setbox\proofabove=\hbox{%
                        \kern-\shortenproofleft\unhbox\proofabove}%
                \shortenproofleft=0pt
        \fi
\fi
\ifdim  \dimen4<0pt
\then   \advance\shortenproofright\dimen4
        \advance\dimen0-\dimen4
        \dimen4=0pt
\fi
\ifdim  \shortenproofright<\wd\proofrulename
\then   \shortenproofright=\wd\proofrulename
\fi
\dimen2=\shortenproofleft \advance\dimen2 by\dimen1
\dimen3=\shortenproofright\advance\dimen3 by\dimen4
\ifproofdots
\then
        \dimen6=\shortenproofleft \advance\dimen6 .5\dimen0
        \setbox1=\vbox to\proofdotseparation{\vss\hbox{$\cdot$}\vss}%
        \setbox0=\hbox{%
                \advance\dimen6-.5\wd1
                \kern\dimen6
                $\vcenter to\proofdotnumber\proofdotseparation
                        {\leaders\box1\vfill}$%
                \unhbox\proofrulename}%
\else   \dimen6=\fontdimen22\the\textfont2 
        \dimen7=\dimen6
        \advance\dimen6by.5\proofrulebreadth
        \advance\dimen7by-.5\proofrulebreadth
        \setbox0=\hbox{%
                \kern\shortenproofleft
                \ifdoubleproof
                \then   \hbox to\dimen0{%
                        $\mathsurround0pt\mathord=\mkern-6mu%
                        \cleaders\hbox{$\mkern-2mu=\mkern-2mu$}\hfill
                        \mkern-6mu\mathord=$}%
                \else   \vrule height\dimen6 depth-\dimen7 width\dimen0
                \fi
                \unhbox\proofrulename}%
        \ht0=\dimen6 \dp0=-\dimen7
\fi
\let\doll\relax
\ifwasinsideprooftree
\then   \let\VBOX\vbox
\else   \ifmmode\else$\let\doll=$\fi
        \let\VBOX\vcenter
\fi
\VBOX   {\baselineskip\proofrulebaseline \lineskip.2ex
        \expandafter\lineskiplimit\ifproofdots0ex\else-0.6ex\fi
        \hbox   spread\dimen5   {\hfi\unhbox\proofabove\hfi}%
        \hbox{\box0}%
        \hbox   {\kern\dimen2 \box\proofbelow}}\doll%
\global\dimen2=\dimen2
\global\dimen3=\dimen3
\egroup 
\ifonleftofproofrule
\then   \shortenproofleft=\dimen2
\fi
\shortenproofright=\dimen3
\onleftofproofrulefalse
\ifinsideprooftree
\then   \hskip.5em plus 1fil \penalty2
\fi
}

\newcommand{\op}{{o}}

\newcommand{\PSbj}{\BBb}
\newcommand{\PObj}{\AAa}

\renewcommand{\to}{\xrightarrow{}}
\newcommand{\ot}{\xleftarrow{}}
\newcommand{\tto}[1]{\xrightarrow{#1}}
\newcommand{\oot}[1]{\xleftarrow{#1}}
\newcommand{\mono}{\rightarrowtail}
\newcommand{\epi}{\twoheadrightarrow}

\newcommand{\Set}{{\sf Set}}

\newcommand{\id}{{\rm id}}

\newcommand{\AAA}{{\cal A}}

\newcommand{\CCC}{{\cal C}}
\newcommand{\DDD}{{\cal D}}
\newcommand{\EEE}{{\cal E}}

\newcommand{\VVV}{{\cal V}}

\newcommand{\XXX}{{\cal X}}
\newcommand{\YYY}{{\cal Y}}

\renewcommand{\Bbb}{\mathbb}
\newcommand{\AAa}{{\Bbb A}}
\newcommand{\BBb}{{\Bbb B}}
\newcommand{\CCc}{{\Bbb C}}
\newcommand{\DDd}{{\Bbb D}}

\newcommand{\FFf}{{\Bbb F}}

\newcommand{\XXx}{{\Bbb X}}

\newcommand{\ZZz}{{\Bbb Z}}

\mathcode`\<="4268 
\mathcode`\>="5269 
\mathchardef\gt="313E 
\mathchardef\lt="313C 

 \def\pushright#1{{
    \parfillskip=0pt            
    \widowpenalty=10000         
    \displaywidowpenalty=10000  
    \finalhyphendemerits=0      
    \leavevmode                 
    \unskip                     
    \nobreak                    
    \hfil                       
    \penalty50                  
    \hskip.2em                  
    \null                       
    \hfill                      
    {#1}                        
    \par}}                      

 \def\qed{\pushright{$\square$}\penalty-700 \smallskip}

\newenvironment{prf}[1]{\begin{trivlist} \item[{\bf ~Proof}#1.]}%
{\qed\end{trivlist}}

\newcommand{\beq}{\begin{equation}}
\newcommand{\eeq}{\end{equation}}
\newcommand{\ba}[1]{\begin{array}{#1}}
\newcommand{\ea}{\end{array}}
\newcommand{\bea}{\begin{eqnarray}}
\newcommand{\eea}{\end{eqnarray}}
\newcommand{\bear}{\begin{eqnarray*}}
\newcommand{\eear}{\end{eqnarray*}}

\newcommand{\Up}[1]{{\Uparrow}{#1}}
\newcommand{\Do}[1]{{\Downarrow}{#1}}
\newcommand{\UD}[1]{{\Updownarrow}{#1}}

\newcommand{\Id}{{\rm Id}}

\newcommand{\lft}[1]{\overleftarrow{#1}}
\newcommand{\rgt}[1]{\overrightarrow{#1}}
\newcommand{\cut}[1]{\overleftrightarrow{#1}} 

\newcommand{\supp}{\mathop{\underrightarrow{\mathrm{lim}}}}
\newcommand{\inff}{\mathop{\underleftarrow{\mathrm{lim}}}}
\newcommand{\limsupp}{\mathop{\overleftarrow{\mathrm{lim}}}}
\newcommand{\liminff}{\mathop{\overrightarrow{\mathrm{lim}}}}

\newcommand{\Dashv}{\mathrel{\mbox{$={\kern-1.1ex}|$}}}
\renewcommand{\vDash}{\mathrel{\mbox{$|{\kern-1.1ex}=$}}}

\newcommand{\Lfun}[1]{{#1}^\#}
\newcommand{\Rfun}[1]{{#1}_\#}
\newcommand{\lan}[1]{{#1}^\ast}
\newcommand{\ran}[1]{{#1}_\ast}
\newcommand{\Lan}[1]{{#1}^\circledast}
\newcommand{\Ran}[1]{{#1}_\circledast}
\newcommand{\Lmon}[1]{\lft{#1}}
\newcommand{\Rmon}[1]{\rgt{#1}}
\newcommand{\LC}{\lan{H}}
\newcommand{\RC}{\ran{H}}
\newcommand{\RLC}{{\Rmon{H}}}
\newcommand{\LRC}{{\Lmon{H}}}

\renewcommand{\paragraph}[1]{\noindent{\bf #1}}

\newcommand{\biarrow}{\looparrowright}

\newcommand{\clsr}[1]{\big[{#1}\big]}
\newtheorem{conjecture}[theorem]{Conjecture}

\DeclareMathOperator{\Img}{Im}
\DeclareMathOperator{\Retr}{Retr}

\newcommand{\place}{\mathord{-}}
\newcommand{\Setaf}[1]{#1\textnormal{-}\Set_{1,\textnormal{free}}}
\newcommand{\Seta}[1]{#1\textnormal{-}\Set}
\newcommand{\abs}[1]{\lvert{#1}\rvert}

\newcommand{\Upper}[1]{{\uparrow}{#1}}
\newcommand{\Lower}[1]{{\downarrow}{#1}}
\newcommand{\UpperLower}[1]{{\updownarrow}{#1}}

\newcommand{\zpz}{{\ZZz_p}}
\newcommand{\od}[2]{\underline{1}{#1}+\underline{p}{#2}}
\newcommand{\odl}[1]{\underline{1}{#1}}
\newcommand{\odr}[1]{\underline{p}{#1}}
\newcommand{\odx}{0}

\graphicspath{{PIC/}}

\title{Towards concept analysis in categories: \\
limit inferior as algebra, 
limit superior as coalgebra
}

\author[1]{Toshiki Kataoka\thanks{Supported by 
Japan Society for the Promotion of Science}}
\author[2]{Dusko Pavlovic\thanks{Supported by AFOSR, AFRL, and UHM.}}
\affil[1]{The University of Tokyo, Tokyo, Japan\\
JSPS Research Fellow\\
  \texttt{toshikik@is.s.u-tokyo.ac.jp}}
\affil[2]{University of Hawaii at Manoa, Honolulu HI, US\\
  \texttt{dusko@hawaii.edu}}
  
  \date{}

\begin{document}

\maketitle

\begin{flushright}
\parbox{10cm}{\small \it It is an open problem whether there exists a sup- and inf-complete category $\AAa''''$ with a sup- and inf-dense embedding $\AAa \to  \AAa''''$ in analogy to the Dedekind completion of an ordered set.}\\[1ex] 
\footnotesize Joachim Lambek \cite[Introduction]{LambekJ:completions}
\\[2ex]
\parbox{10cm}{\small \it No Lambek extension of the one-object category $\ZZz_4$ has finite limits.}
\\
{\footnotesize John Isbell \cite[Thm. 3.1]{IsbellJ:no-lambek}}
\end{flushright}

\begin{abstract}
While computer programs and logical theories begin by declaring the concepts of interest, be it as data types or as predicates, network computation does not allow such global declarations, and requires \emph{concept mining}\/ and \emph{concept analysis}\/ to extract shared semantics for different network nodes. Powerful semantic analysis systems have been the drivers of nearly all paradigm shifts on the web. In categorical terms, most of them can be described as bicompletions of enriched matrices, generalizing the Dedekind-MacNeille-style completions from posets to suitably enriched categories. Yet it has been well known for more than 40 years that ordinary categories themselves in general do not permit such completions. Armed with this new semantical view of Dedekind-MacNeille completions, and of  matrix bicompletions, we take another look at this ancient mystery. It turns out that simple categorical versions of the \emph{limit superior}\/ and \emph{limit inferior}\/ operations characterize a general notion of Dedekind-MacNeille completion, that seems to be appropriate for ordinary categories, and boils down to the more familiar enriched versions when the limits inferior and superior coincide. This explains away the apparent gap among the completions of ordinary categories, and broadens the path towards categorical concept mining and analysis, opened in previous work. 
\end{abstract}

\section{Introduction}

\subsection{Problem of concept mining and analysis}
\begin{figure}[htbp]
\centering
\includegraphics[height=3.8cm]{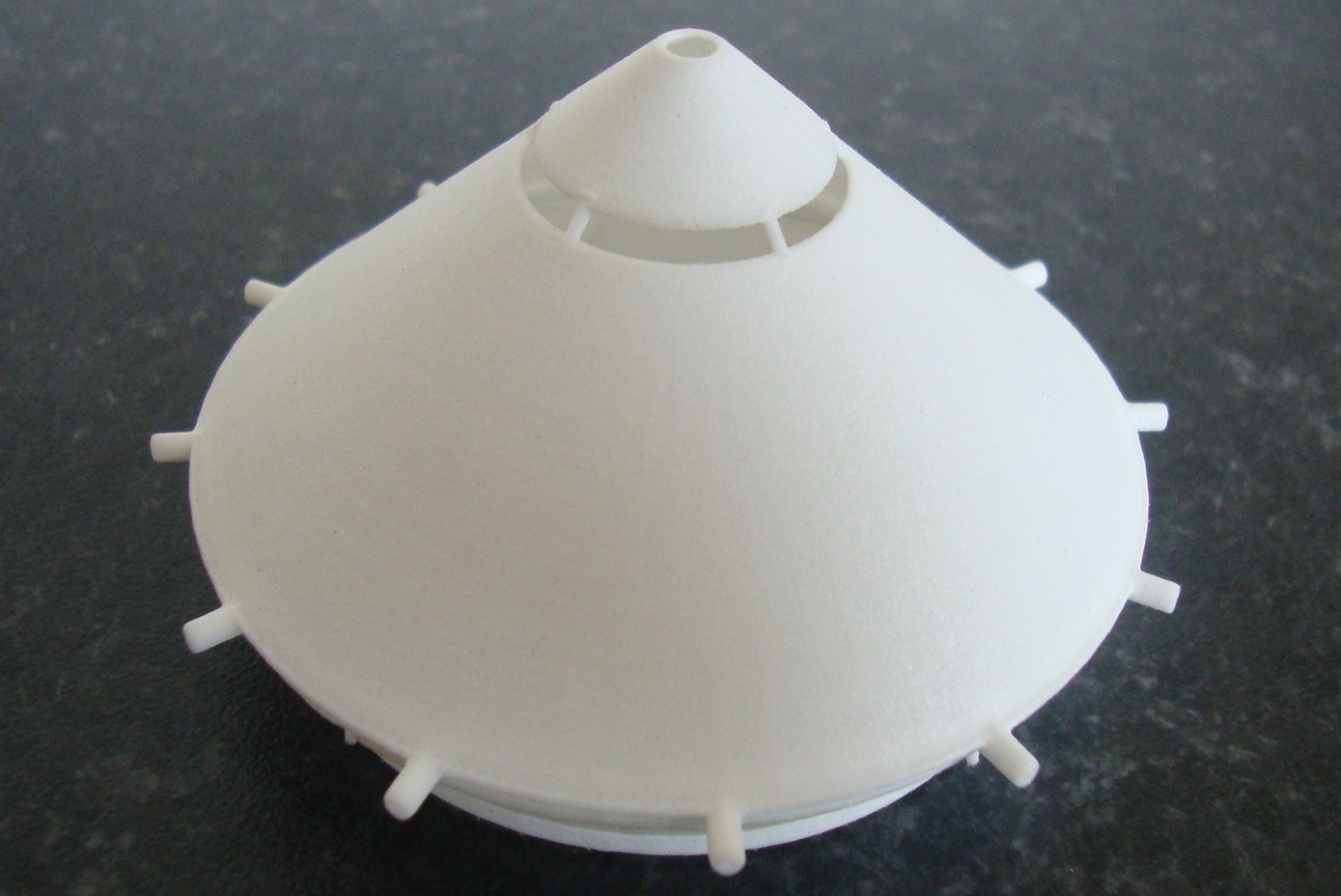}\ \ \ \includegraphics[height=3.8cm]{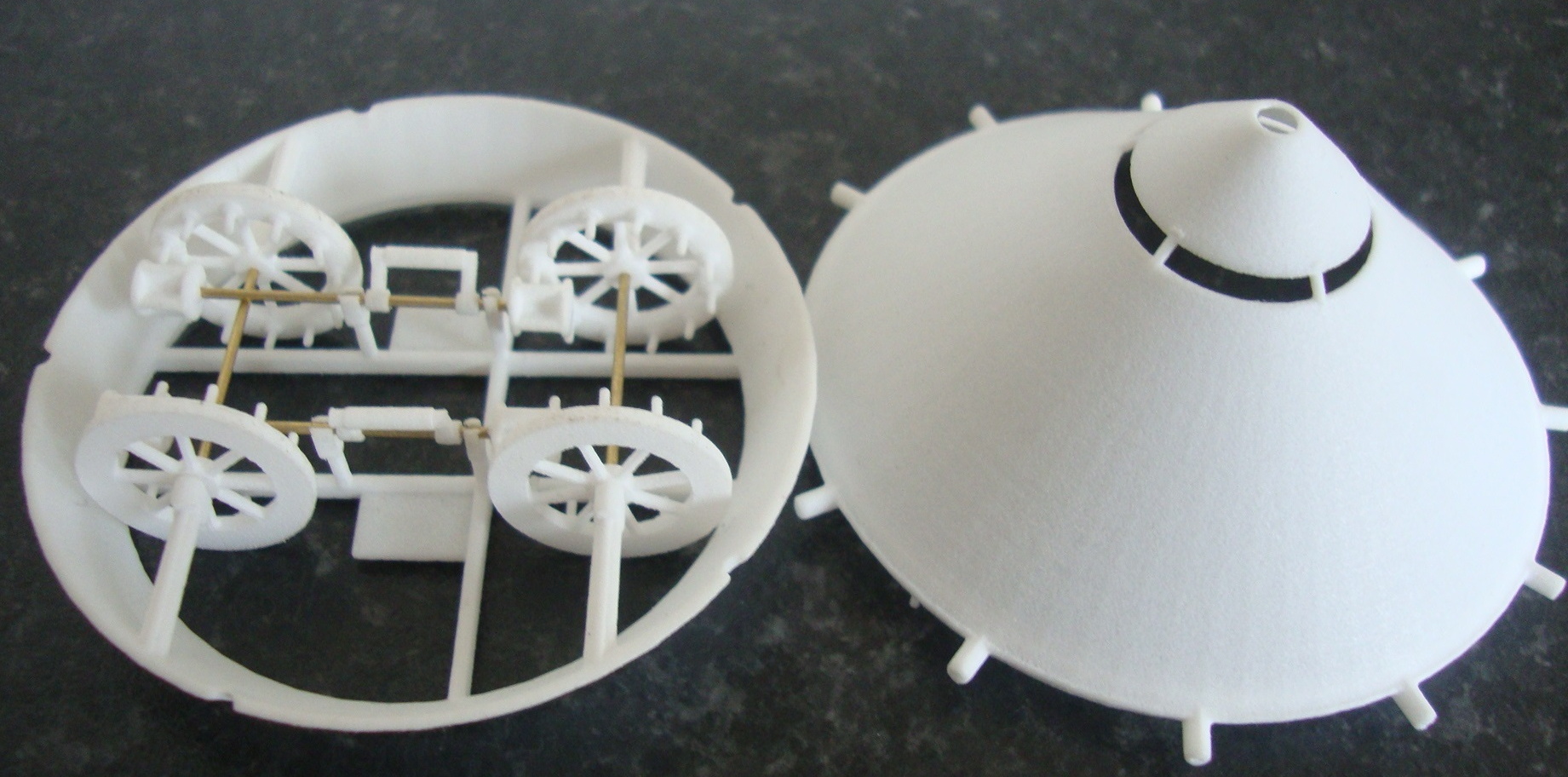}
\caption{Unidentified object: The external and the internal view}
\label{figone}
\end{figure}

Suppose you come across upon  the object depicted  in Fig.~\ref{figone}. The conic top is easily removed to uncover the mechanism on the right. 
What is this thing?

You would surely approach the problem from both directions at once: on one hand, you would look how the parts  fit together and try to discern the \emph{structural components}\/ of the device; on the other hand, you would twiddle with some parts and watch what moves together, trying to figure out the \emph{functional modules}. The parts that move together may not be next to each other, but they probably belong to the same functional module. The parts that are related structurally are more likely to be related functionally. If you manage to discern some distinct components corresponding to distinct functionalities, then each such component-function pair will presumably correspond to a \emph{concept}\/ conceived by the designer of the device. By analyzing the device you will extract the designer's idea.

Similar analyses are formalized under different names in different research communities: some speak of \emph{concept analysis}, some of \emph{knowledge acquisition}, \emph{semantic indexing}, or \emph{data mining} \cite{Carpineto-Romano:book,GardenforsP:conceptual,data-mining-handbook}. The application domains and the formalisms vary very widely, from mathematical taxonomy \cite{Jardine-Sibson}, through text analysis \cite{text-mining} and pattern recognition \cite{Duda:pattern}, to web search and recommender systems \cite{recommender-handbook}. The importance of formalizing and implementing concept analysis grew rapidly with the advent of the web, as almost anything found on the web requires some sort of concept mining and analysis, not only because there are no global semantical declarations, and the meaning has to be extracted from the network structure \cite{PavlovicD:CSR08}, but also to establish trust \cite{PavlovicD:FAST10}. Diverse toy examples of such concept analysis tasks, motivating the modeling approach extended in this paper, can be found in \cite{PavlovicD:QI08,PavlovicD:FAST10,PavlovicD:ICFCA12,PavlovicD:Samson13}.

The analytic process that a formal concept analyst may initiate upon  an encounter with the unidentified object from Fig.~\ref{figone} is thus not all that different from what a curious child would do: they would both start by recording the observed components on one hand, and the observed functionalities on the other, and they would note which components are related to which functionalities. With the 'yes-no' relations, the formal version of this process leads to the simple and influential method that goes under the name \emph{Formal Concept Analysis  (FCA)} \cite{FCA-book,FCA-foundations}. If the relations between the components and the functionalities are quantified by real numbers and stored in \emph{pattern matrices}, then the analysis usually proceeds by the methods of statistics and linear algebra, and goes under the name \emph{Principal Component Analysis (PCA)} \cite{Jolliffe:PCA}, or \emph{Latent Semantic Analysis (LSA)} \cite{LSA:book}, etc. It performs the singular value decomposition of the pattern matrix, and thus mines the concepts as the eigenspaces of the induced linear operators. 

Interestingly, if you wanted to record that the unidentified device has 4 identical wheels, and that each wheel has 12 identical cogs, and that two of the wheels are related to two different functionalities, driving and steering, you would be led beyond the familiar concept mining approaches. While the experts in these approaches would surely figure out multiple tricks to record what is needed (e.g. by using multi-level pattern matrices), the straightforward approach leads beyond the FCA matrices of 0s and 1s, and beyond the LSA matrices of real numbers, to matrices of sets between components on one hand, matrices of sets between functionalities on the other hand, and matrices of sets between components and functionalities in-between. You would construct a category of components, a category of functionalities, and a profunctor/distributor between them. If the cog is recognized as a part, then a coproduct of 12 cogs would be embedded in each wheel. If the cogs are attached with rivets, then their morphisms may not be monic, since the distinctions of some of their parts may be obliterated through deformations. So why have such categorical models not been used in concept analysis?

%

Many of the concept mining approaches derived from LSA are instances of \emph{spectral decomposition}\/ \cite{Azar}. Formalized in terms of \emph{enriched category theory}\/ \cite{KellyGM:enriched-book}, the problem of concept mining turns out to be an instance of a general spectral decomposition problem \cite{PavlovicD:FAST10,PavlovicD:ICFCA12}, which can also be viewed as a problem of \emph{minimal bicompletion}\/ of a suitably enriched matrix \cite{PavlovicD:Samson13}. Even the standard linear algebra of LSA seems to be an instance of such bicompletion, over a suitable category\footnote{not poset!} of real numbers. The problem of minimal bicompletions of enriched categories, which subsume the Dedekind-MacNeille completions of posets, is the special case, arising when a category itself is viewed as a matrix. Instantiated to categories enriched over sets, also known as \emph{'ordinary'}\/ categories, this turned out to be a strange problem, as suggested by the quotations at the very beginning of the paper. Maybe this is the reason for the notable absence of ordinary categories in the extensive concept mining toolkits? We sketch the problem of bicompletions of ordinary categories  of in the next section.

\subsection{Problem of minimal bicompletions of matrices and categories}\label{Sec-problem}
Throughout the paper, we assume familiarity with the basic concepts of category theory, e.g. at the level of \cite{CWM}. To understand the general approach to concept mining through minimal bicompletions, explained in this section, the reader may need some ideas about enriched categories as well, e.g. as presented in \cite{KellyGM:enriched-book}. Beyond this section, the rest of the paper will be about ordinary categories.

Suppose that we have thus proceeded as in the preceding section, and built a category of components $\PObj$ and a category of functionalities $\PSbj$. If we have recorded just the inclusion relations, then each of these categories is a poset, i.e. enriched over the ordered monoid $(\{0,1\}, \wedge, 1)$. If we have recorded the distances among the components on one hand, and among the functionalities on the other, then our categories are metric spaces \cite{LawvereFW:metric}, viewed as categories enriched over the monoidal poset $
\left([0,\infty], +, 0\right)$. If we capture the components and the functionalities as ordinary categories, then $\PObj$ and $\PSbj$  are enriched over the monoidal category $(\Set, \times, 1)$. 

\subsubsection{The setting of minimal bicompletion}
The relationships between the components and the functionalities will be expressed as a $\VVV$-enriched functor $\Phi: \PObj^\op \times\PSbj \to \VVV$, where $\VVV$ is the enriching category, such as $\{0,1\}, [0,\infty]$ or $\Set$ above. We call such $\VVV$-enriched functor a \emph{matrix}. In particular, given a $\VVV$-matrix $\Phi: \PObj^\op\times \PSbj \to \VVV$ we derive its extensions as in Fig.~\ref{fig-deriv}.
\begin{figure}[htbp]
\centering
$
\prooftree
 \prooftree
 \prooftree
 \Phi:\PObj^\op \times \PSbj \to \VVV
 \justifies
\Lfun \Phi : \PObj \to \left(\VVV^\PSbj\right)^\op\qquad \qquad
\Rfun \Phi : \PSbj \to \VVV^{\PObj^\op}  
 \endprooftree
 \justifies
\lan \Phi : \VVV^{\PObj^\op} \to \left(\VVV^\PSbj\right)^\op
 \qquad \qquad \ran \Phi : \left(\VVV^\PSbj\right)^\op \to \VVV^{\PObj^\op}  \endprooftree
 \justifies
 \lft \Phi = \ran \Phi \lan \Phi : \VVV^{\PObj^\op} \to \VVV^{\PObj^\op}  \qquad \qquad   \rgt \Phi = \lan \Phi \ran \Phi : \left(\VVV^{\PSbj}\right)^\op \to \left(\VVV^{\PSbj}\right)^\op
 \endprooftree  
 $
\caption{Deriving the two extensions and the two kernels of a matrix $\Phi$}
\label{fig-deriv}
\end{figure}
The functors $\Lfun \Phi$ and $\Rfun \Phi$ are the transpositions of $\Phi$. The presheaves in the form $\Rfun \Phi b$ and the postsheaves in the form $\Lfun \Phi a$ are called \emph{$\Phi$-representable}. The functors $\lan \Phi$ and $\ran \Phi$ are the Kan extensions \cite[Ch.~4]{KellyGM:enriched-book} of $\Lfun \Phi$ and $\Rfun \Phi$. Since they form an adjunction, their composite $\lft \Phi$ is a monad and $\rgt \Phi$ is a comonad.

\begin{figure}[htbp]
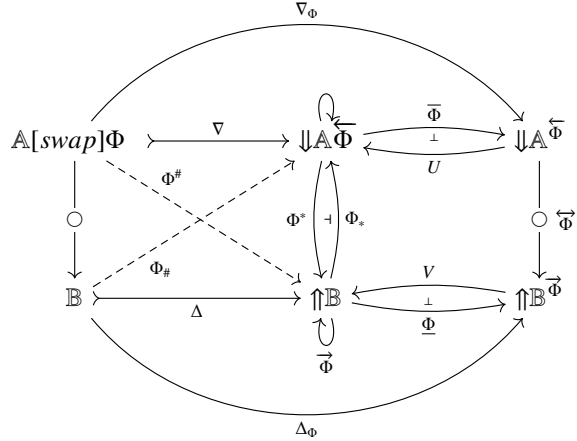

\centering
\begin{tikzar}{}
\PObj \arrow[tail]{rr}{\nabla} \arrow[bend left=50]{rrrr}{\nabla_\Phi}  \arrow{dd}[description]{\bigcirc}
[swap]{\Phi\ \ }
\arrow[dashed]{ddrr}[near start]{\Lfun \Phi} \&\& \Do \PObj  \arrow[bend right=13]{dd}[swap]{\lan \Phi} \arrow[loop above]{} {\Lmon \Phi} \arrow[bend left=10]{rr}{\overline\Phi}[swap]{\scriptscriptstyle\bot} \&\& \Do \PObj^{\Lmon \Phi} \arrow[bend left=10]{ll}{U} 
\arrow{dd}{\ \ \cut\Phi}[description]{\bigcirc}
\\
\\
\PSbj \arrow[tail]{rr}[swap]{\Delta} \arrow[bend right=50]{rrrr}[swap]{\Delta_\Phi} \arrow[dashed]{uurr}[near start,swap]{\Rfun\Phi} \& \& \Up \PSbj \arrow[bend right=13]{uu}{\dashv}[swap]{\ran\Phi} \arrow[loop below]{}{\Rmon\Phi} \arrow[bend right=10]{rr}{\scriptscriptstyle\bot}[swap]{\underline \Phi} \&\& \Up \PSbj^{\Rmon \Phi} \arrow[bend right=10]{ll}[swap]{V}
\end{tikzar}
\caption{Minimal bicompletion of a matrix $\Phi$}
\label{fig-extensions}
\end{figure}

When the enrichment is clear from the context, it is convenient to abbreviate the matrix $\PObj^\op\times \PSbj \to \VVV$ to $\PObj \biarrow \PSbj$ and the completions $\VVV^{\PObj^\op}$ and $\left(\VVV^\PSbj\right)^\op$ to $\Do \PObj$ and $\Up \PSbj$ respectively, so that the derivations in Fig.~\ref{fig-deriv} give the diagram in Fig.~\ref{fig-extensions}
where $\nabla$ and $\Delta$ are the Yoneda embeddings \cite[Sec.~2.4]{KellyGM:enriched-book}. The monad $\Lmon \Phi =  \ran \Phi \lan \Phi$, induced by the  Kan extensions $\lan \Phi \dashv \ran \Phi: \Up \PSbj \to \Do \PObj$, induces the category of (Eilenberg-Moore) algebras $\Do \PObj^{\lft \Phi}$, whereas the comonad $\Rmon \Phi  =   \lan \Phi \ran \Phi$ induces $\Up \PSbj^{\rgt \Phi}$. The functor $\nabla_\Phi = \overline \Phi \circ \nabla$ maps $\PObj$ to free $\lft \Phi$-algebras generated by the representable presheaves,
whereas the functor $\Delta_\Phi = \underline \Phi \circ \Delta$ maps $\PSbj$ to cofree $\rgt \Phi$-coalgebras cogenerated by the representable postsheaves.

\subsubsection{Familiar cases}
When $\VVV = \{0,1\}$, the $\VVV$-enriched categories $\PObj$ and $\PSbj$ are posets. Then $\Do \PObj$ consists of antitone maps $\lft L:\PObj^\op \to \{0,1\}$, or equivalently of the lower-closed sets in $\PObj$, whereas $\Up \PSbj$ consists of the monotone maps $\rgt U: \PSbj \to \{0,1\}$, or equivalently of the upper-closed sets in $\PSbj$. The Yoneda embedding $\nabla: \PObj \to \Do \PObj$ is then the supremum (or join) completion, and the $ \Delta:\PSbj \to \Up \PSbj$ is the infimum (or meet) completion. A matrix $\Phi: \PObj^\op \times \PSbj \to \{0,1\}$ corresponds to a subset of the product poset which is lower closed in $\PObj$ and upper closed in $\PSbj$. Its extensions are then
\bea\label{eq-poset-extensionone}
\lan \Phi \lft L & = & \left\{u\in \PSbj\ |\ \forall x.\ \lft L(x) \Rightarrow \Phi(x,u)\right\}\\  
\ran \Phi \rgt U & = & \left\{\ell\in \PObj\ |\ \forall y.\  \rgt U(y)\Rightarrow \Phi(\ell, y)\right\} \label{eq-poset-extensiontwo}
\eea
Intuitively, $\lan \Phi\lft L$ can be construed as the set of upper bounds in $\PSbj$ of the $\Phi$-image of the lower set $\lft L$, whereas $\ran \Phi \rgt U$ can be construed as the set of $\Phi$-lower bounds of the upper set $\rgt U$. The operator $\Lmon \Phi = \ran \Phi \lan \Phi$ thus maps each lower set $\lft L$ to the set of the $\Phi$-lower bounds of the set of its $\Phi$-upper bounds; whereas the operator $\Rmon \Phi = \lan \Phi \ran \Phi$ maps each upper set $\rgt U$ to the set of the $\Phi$-upper bounds of its $\Phi$-lower bounds. Both operators are thus \emph{closure operators}. Their lattices of closed sets $\Do \PObj^{\Lmon \Phi}$ and $\Up \PSbj^{\Rmon \Phi}$ turn out to be isomorphic, and form the \emph{nucleus} of $\Phi$ \cite{PavlovicD:ICFCA12}. A $\Lmon \Phi$-closed set in $\PObj$ and the corresponding $\Rmon \Phi$-closed set in $\PSbj$, of course, completely determine each other, but the most informative presentation carries both, as Dedekind-style cuts.  When $\PObj = \PSbj$  and $\Phi\subseteq \PObj^\op \times \PObj$ is the partial ordering 
\bea\label{eq-order}
\Phi(x,y) & \iff & x\leq y
\eea
then the nucleus is just the Dedekind-MacNeille completion $\UD \PObj$ of the poset $\PObj$ \cite{MacNeille}. This is the \emph{minimal}\/ bicompletion, in the sense that the embedding $\PObj \to \UD\PObj$ preserves any suprema and infima that $\PObj$ may already have, and only adds those that do not yet exist \cite[III.3.11]{MacNeille,Banaschewski-Bruns,JohnstoneP:stone}. The consequence of this minimality is that every element of the completion $\UD\PObj$ is \emph{both}\/ a supremum and an infimum of the elements of $\PObj$. The nucleus of a $\{0,1\}$-matrix is a minimal bicompletion in a similar sense, as are the nuclei of $[0,1]$-matrices, and of $[0,\infty]$-matrices\footnote{Since the monoidal posets $\big([0,\infty], +, 0\big)$ and $\big([0,1], \times , 1\big)$ are isomorphic as monoidal categories, all statements about categories enriched over them transfer trivially. However, isomorphisms are not always trivial phenomena. E.g., the Laplace transform is an isomorphism, which maps differential operations into algebraic operations, and thus allows solving differential equations as algebraic equations, and mapping back the solutions \cite{PavlovicD:LICS98}. In a similar way, it often happens that a distance space presentation of a data pattern, enriched over $\big([0,\infty], +, 0\big)$, displays some geometric content, whereas an isomorphic proximity lattice presentation of the same data pattern, enriched over $\big([0,1], \times , 1\big)$, displays some generalized order structure, not apparent in the first interpretation.} $
\big([0,1], \times , 1\big)$: the nuclei give the \emph{semantic bicompletions}\/ of matrices, uncovering their concepts \cite{PavlovicD:ICFCA12,PavlovicD:Samson13}.

\subsubsection{The trouble with ordinary categories}
Our main concern in the present paper are the minimal bicompletions of matrices and categories enriched over $(\Set,\times, 1)$. Categories enriched in $\Set$ are usually called \emph{ordinary}\/ categories. $\Set$-matrices are variably called \emph{profunctors}\/ or \emph{distributors}. We increase  the wealth of terminology by calling them \emph{matrices}. The functors $\lft \alpha \in \Do \AAa = \Set^{\PObj^\op}$ are called \emph{presheaves}. The functors $\rgt \beta \in \Up \BBb = \left(\Set^{\PSbj}\right)^\op$ are usually called \emph{covariant functors to $\Set$}, but we call them \emph{postsheaves}. We use without further explanation the well known fact \cite{SGA1,MacLane-Moerdijk} that presheaves are equivalent to discrete fibrations, and that postsheaves are equivalent with discrete opfibrations.

We also call the categorical limits the \emph{infima}, and the categorical colimits the \emph{suprema}, following Lambek's 1966 Lectures on Completions of Categories \cite{LambekJ:completions}, quoted at the beginning of this paper. The Yoneda embeddings $\nabla: \PObj \to \Do \PObj$ and $ \Delta:\PSbj \to \Up \PSbj$ are then again, respectively, the supremum and the infimum completion, this time of the categories $\PObj$ and $\PSbj$. The transposes $\Lfun \Phi$ and $\Rfun \Phi$ now extend to the adjunction $\lan \Phi \dashv \ran \Phi: \Up \PSbj \to \Do \PObj$, which are defined similarly to (\ref{eq-poset-extensionone}--\ref{eq-poset-extensiontwo}). More precisely, the mappings between the $\PObj$-presheaves and $\PSbj$-postsheaves 
\[\hspace{8em}
\prooftree
\lft \alpha : \PObj^\op \to \Set 
\justifies
\lan \Phi \lft \alpha : \PSbj \to \Set
\endprooftree
\hspace{8em}
\prooftree
\rgt \beta : \PSbj \to \Set 
\justifies
\ran \Phi \rgt \beta : \PObj^\op \to \Set
\endprooftree
\]
are defined as follows
\bea\label{eq-cat-extensionone}
\lan \Phi \lft \alpha\left(u\right) & = &\inff_{x\in \PObj} \Big(\lft \alpha(x) \Rightarrow \Phi(x,u)\Big)\  \ \ =\ \ \  \Do \PObj \left(\lft \alpha, \Lfun{\Phi} u\right)\\  
\ran \Phi \rgt \beta\left(\ell\right) & = & \inff_{y\in \PSbj} \Big(\rgt \beta(y)\Rightarrow \Phi(\ell, y)\Big) \ \ \  = \ \ \ \Up \PSbj  \left(\Rfun{\Phi} \ell, \rgt \beta\right) \label{eq-cat-extensiontwo}
\eea
Here we write $X\Rightarrow Y$ for the set exponents $Y^X$ not only because the multiple exponents tend to 'fly away' in the latter notation, but also to emphasize the parallel with (\ref{eq-poset-extensionone}--\ref{eq-poset-extensiontwo}). When $\PObj = \PSbj$ is the same category, and $\Phi = H : \PObj^\op \times \PObj \to \Set$ is the hom-set matrix, then $\lan H\lft \alpha (u)$ is the set of (right) cones from the presheaf $\lft \alpha$, viewed as a diagram, to the object $u$ as the tip of the cone. Dually, $\ran H\rgt \beta (\ell)$ is the set of (left) cones from the tip $\ell$ to the diagram $\rgt \beta$. For a general matrix $\Phi$, thinking of the elements of each set $\Phi(a,b)$ as 'arrows' from  $a\in \PObj$ to $b\in \PSbj$ also allows thinking of $\rgt \varrho \in \lan \Phi\lft \alpha (u)$ as a (right) 'cone' from the diagram $\lft \alpha$ in $\PObj$ to the tip $u\in \PSbj$, and of $\lft \lambda \in \ran \Phi\lft \beta (\ell)$ as a (left) 'cone' from the tip $\ell\in \PObj$ to a diagram $\lft \beta$ in $\PSbj$. The presheaves and postsheaves of \eqref{eq-cat-extensionone} and \eqref{eq-cat-extensiontwo} thus generalize the lower and the upper sets of \eqref{eq-poset-extensionone} and \eqref{eq-poset-extensiontwo}. 

At the very beginning of his lectures, Lambek raised the question of the Dedekind-MacNeille completion of a category, and left it open. He did not raise the general question of semantic completions of matrices (profunctors, or distributors) only because the semantical impact was not clear at the time; but the general situation from Fig.~\ref{fig-extensions} was well known. Lambek's open question of the Dedekind-MacNeille completion of a category was closed by Isbell a couple of years later, who showed in \cite[Sec.~3]{IsbellJ:no-lambek} that already the group $\ZZz_4$, viewed as a category with a single object, cannot have a completion generated both by the suprema and by the infima.

However, taking a broader semantical view, and seeking semantic completions of matrices, shows that the story does not really end with Isbell's counterexample. A semantic completion of a matrix, relating, say, the parts and the moves observed within a device like the one on Fig.~\ref{figone}, should uncover the concepts underlying the design of the device. These concepts are expressed through the structural component of the device, and through its functional units. When the matrix is enriched over a monoidal poset, then there is a one-to-one correspondence between the structural components and the functional modules, and they form the nucleus of the matrix \cite{PavlovicD:ICFCA12,PavlovicD:Samson13}. In reality, though,  a single structural component may play a role in several functional modules, and vice versa. While the posetal enrichment cannot capture this, the enrichment in sets, or in a proper category of real numbers, can record how many copies of a given a part are used for a certain function. Modeled in this way, the spaces of structural components  and of functional modules will not be isomorphic. The concepts will not be  uncovered as a single category of component-function pairs, like in the posetal case, but as a nontrivial matrix relating some component-concepts approximated by their functionalities with some function-concepts approximated by the components that perform them. 

\subsubsection*{Contributions}
To spell this out, we consider the following technical questions:
\begin{enumerate}[(a)]
\item What kind of completions of a given matrix $\Phi: \PObj \biarrow \PSbj$ are provided by the categories $\Do \PObj^{\lft \Phi}$ and $\Up \PSbj^{\rgt \Phi}$? (The idea is that the former captures  the component-concepts, the latter the function-concepts.) 
\item What kind of matrix $\cut \Phi: \Do\PObj^{\lft \Phi} \biarrow \Up\PSbj^{\rgt \Phi}$ is the minimal bicompletion of $\Phi: \PObj \biarrow \PSbj$? (Capturing the relations between the component-concepts and the function-concepts.)
\end{enumerate}
Our approach to these questions is based on a new family of limits and colimits, introduced in the next section. It seems intuitive and appropriate to call them \emph{limit inferior}, and \emph{limit superior}. For consistency, we also revert, albeit just for the duration of this paper\footnote{We hope that our terminological contributions, advancing from 'profunctors' and 'distributors' to 'matrices' and from 'covariant functors to $\Set$' to 'postsheaves', as well as retreating from 'limits' to 'infima' and from 'colimits' to 'suprema', will not end up being the central features of the paper.}, from \emph{limits}\/ and \emph{colimits}\/ to \emph{infima}\/ and \emph{suprema}, following Lambek \cite{LambekJ:completions}. The reader is reminded that in in posets 
\begin{itemize}
\item the limit inferior is the \emph{supremum}\/ of the lower bounds of a set, whereas
\item the limit superior is the \emph{infimum}\/ of the upper bounds.
\end{itemize}
Mutatis mutandis, the categorical concepts will behave similarly.

\subsubsection*{Overview of the paper}
In Sec.~\ref{Sec-two}, we propose the answers to the above question. Sec.~\ref{Sec-two-zero} spells out the preliminaries. Sec.~\ref{Sec-two-one} defines categorical limits inferior and superior and characterizes their completions. Sec.~\ref{Sec-two-two} proposes an answer to question (a) above. Sec.~\ref{Sec-two-three} proposes an answer to question (b) above. In Sec.~\ref{Sec-three} we study some simple examples, illustrating and validating the introduced concepts. Sec.~\ref{Sec-three-one} describes a monadicity workflow useful for analyzing the examples. Sec.~\ref{Sec-three-two} characterizes completions of constant matrices. Sections~\ref{Sec-three-three}  and \ref{Sec-three-four} characterize completions of the matrices representing groups or posets, respectively. Sec.~\ref{Sec-three-five}  characterizes completions of a vector in the group $\ZZz_p$ for a prime $p$. Sec.~\ref{Sec-four} closes the paper, to some extent.

Due to the space constraints of this conference paper and the scope of the presented material, all proofs and many lemmas had to be moved into the appendices. Full details will require a significantly longer paper.

\section{Categorical limit inferior and limit superior}\label{Sec-two}

\subsection{Preliminaries}\label{Sec-two-zero}
Although suprema and infima are very basic concepts, familiar to most readers, and easily found in \cite[Sec.~III.3--4]{CWM}, we spell them out here not only to introduce the notation and practice using the words \emph{infimum}\/ and \emph{supremum}\/ instead of limit and colimit, but also to align these familiar definitions with the variations needed to define the \emph{limit superior}\/ and the \emph{limit inferior}.

Let $\CCc$ and $\DDd$ be categories and $\CCc^\DDd$ the category of functors between them, with natural transformations as morphisms. Let $\Box: \CCc \to \CCc^\DDd$ be the functor taking each object of $x$ of $\CCc$ to the constant functor $\Box x:\DDd\to \CCc$, which maps all objects of $\DDd$ to $x\in \CCc$ and all morphisms of $\DDd$ to $\id_x$. 

The suprema and the infima in $\CCc$ can be defined as, respectively, the left and the right adjoint of the constant functor, i.e.
\bear
\supp \dashv \Box \dashv \inff & : & \CCc^\DDd\to \CCc
\eear
These adjunctions can be viewed as the natural bijections
\bea
\CCc^\DDd(F, \Box x) & \cong & \CCc\left(\supp F, x\right)\label{def-sup}\\
\CCc^\DDd(\Box x, F) & \cong & \CCc\left(x, \inff F\right)\label{def-inf}
\eea

It is well known that the Yoneda embeddings realize the $\supp$ and $\inff$-completions \cite[Sec.~X.6]{CWM}:
\begin{itemize}
\item $\nabla : \CCc \to \Do \CCc$ is the $\supp$-completion of $\CCc$, whereas
\item $\Delta : \CCc \to \Up \CCc$ is the $\inff$-completion of $\CCc$
\end{itemize}
where 
\begin{itemize}
\item $\Do \CCc$ denotes the category $\Set^{\CCc^\op}$ of $\CCc$-\emph{presheaves}, or equivalently\footnote{The equivalence between the "indexed" and "fibered" versions of sheaves lies at the heart of Grothendieck's descent theory  \cite[VI]{SGA1}, but also generalizes to substantially different purposes \cite{PavlovicD:CTCS97,PavlovicD:QPL09}.} the category of discrete fibrations over $\CCc$,
\item $\Up \CCc$ denotes the category $\left(\Set^\CCc\right)^\op$ of $\CCc$-\emph{postsheaves}, or equivalently the opposite category of discrete opfibrations over $\CCc$.
\end{itemize}

For completeness, we note the following well known and routinely checkable fact.
\begin{lemma} \label{lemma-completion-algebra}
Given a functor $F:\DDd\to \CCc$, consider the presheaf and the postsheaf
\[
\big(\lft F :\ \CCc/\!\!/F \to \CCc\big) \ \in \Do \CCc\qquad \qquad \qquad
\big(\rgt F :\ F/\!\!/\CCc \to \CCc\big) \ \in \Up \CCc 
\]
where $\CCc/\!\!/F$ is the category of connected components of the comma category $\CCc/F$ from $\Id_\CCc$ to $F$, whereas $F/\!\!/\CCc$ is the category of connected components of the comma category $F/\CCc$ the other way around \cite[Sections~II.6 and IX.3]{CWM}. Then
\[
\supp F =  \supp \lft F\qquad \qquad \qquad
\inff F  =  \inff \rgt F
\]
\end{lemma}

\noindent \textbf{Notations} have been introduced in Sec.~\ref{Sec-problem}, especially in Figures~\ref{fig-deriv} and \ref{fig-extensions}. The next section considers the special case $\Phi = H:\CCc \biarrow \CCc$ of the matrix of hom-sets of a category.

\subsection{Limit inferior and limit superior over a category}\label{Sec-two-one}
\begin{definition}
For arbitrary categories $\CCc$ and $\DDd$ we define 
\begin{itemize}
\item the category of \emph{left saturated diagrams} $\CCc^\DDd_\Downarrow$ to consist of
\begin{itemize}
\item objects $\lvert \CCc^\DDd_\Downarrow \rvert = \lvert \CCc^\DDd \rvert$
\item morphisms $\CCc^\DDd_\Downarrow\left( F, G \right) = \Do \CCc\left( \RC \rgt F , \RC \rgt G\right)$
\end{itemize}
\item the category of \emph{right saturated diagrams} $\CCc^\DDd_\Uparrow$ to consist of
\begin{itemize}
\item objects $\lvert \CCc^\DDd_\Uparrow \rvert = \lvert \CCc^\DDd \rvert$
\item morphisms $\CCc^\DDd_\Uparrow\left( F, G \right) = \Up \CCc\left( \LC \lft F , \LC \lft G \right)$
\end{itemize}
\end{itemize}
\end{definition}

\begin{definition}
In a category $\CCc$ we define 
\begin{itemize}
\item the \emph{limit inferior}\/ operation $\liminff$ over left diagrams from $\DDd$ by the adjunction
\bear
\liminff \dashv \Box & : & \CCc \to \CCc^\DDd_\Downarrow
\eear
which can be viewed as the natural bijection
\bear
\CCc^\DDd_\Downarrow \left(F, \Box x\right) & \cong & \CCc\left( \liminff F, x\right)
\eear
\item the \emph{limit superior}\/ operation $\limsupp$  over right diagrams from $\DDd$ by the adjunction
\bear
\Box \dashv \limsupp & : & \CCc^\DDd_\Uparrow \to \CCc
\eear
which can be viewed as the natural bijection
\bear
\CCc^\DDd_\Uparrow \left(\Box x, F\right) & \cong & \CCc\left(x, \limsupp F\right)
\eear
\end{itemize}
\end{definition}

\paragraph{Remarks.} Note that the operations $\liminff$ and $\limsupp$ are defined over \emph{arbitrary}\/ diagrams. Indeed, the objects of the categories of \emph{saturated}\/ diagrams are arbitrary diagrams; the saturation is imposed on them in the definitions of the morphisms in these categories. 

The ordinary infima and suprema are also defined over arbitrary diagrams, but differently: a supremum of a diagram is equal to the supremum of the  induced presheaf; and the infimum of a diagram is equal to the infimum of the induced postsheaf, as stated in  Lemma \ref{lemma-completion-algebra}. This is analogous to lattices, where a supremum of a set is equal to the supremum of its lower closure, whereas the infimum of a set is the infimum of the upper closure. However,  the limit inferior of a diagram is the supremum of the presheaf induced by the postsheaf induced by the diagram; and the limit superior is the infimum of the postsheaf induced by the presheaf induced by the diagram. In a partially ordered set, the limit inferior of a set is the join of the lower bounds of all of its upper bounds; whereas the limit superior of a set is the meet of the upper bounds of all of its lower bounds.

\begin{lemma}\label{lemma-repres}
Every representable presheaf $\nabla x$ is a free algebra in  $\Do \CCc^{\LRC }$, with $\nabla x\stackrel\eta \cong\lft H \nabla x$. Every representable postsheaf $\Delta x$ is a cofree coalgebra in $\Up \CCc^{\RLC}$, with $\rgt H \Delta x\stackrel \varepsilon \cong \Delta x$.
\end{lemma}

\begin{proposition}\label{prop-represgen}
Every $\LRC$-algebra is a limit inferior in $\Do \CCc^{\LRC }$ of representable  presheaves, viewed as $\LRC$-algebras. Every $\RLC$-coalgebra is a limit superior in $\Up \CCc^{\RLC}$ of representable postsheaves, viewed as $\RLC$-coalgebras. 
\end{proposition}

\begin{corollary}\label{cor-complete}
$\Do \CCc^{\LRC }$ is $\liminff$-complete. $\Up \CCc^{\RLC}$ is $\limsupp$-complete.  
\end{corollary}

\begin{theorem}\label{thm-completion}
The extended Yoneda embeddings  realize the limit inferior and limit superior completions:
\begin{itemize}
\item $\nabla_H : \CCc \tto\nabla \Do \CCc\tto{\overline H} \Do \CCc^{\lft H}$ is the $\liminff$-completion of $\CCc$, whereas
\item $\Delta_H : \CCc \tto{\Delta} \Up \CCc\tto{\underline H} \Up \CCc^{\rgt H}$ is the $\limsupp$-completion of $\CCc$.
\end{itemize}
\end{theorem}
\subsection{Limit inferior and limit superior over a matrix}\label{Sec-two-two}
Given a category $\CCc$, Lemma~\ref{lemma-completion-algebra} implies that the suprema and the infima, defined by \eqref{def-sup} and \eqref{def-inf} respectively, can be viewed as the left and the right adjoint of the corresponding Yoneda embeddings:
\begin{center}
\begin{tikzar}{}
\CCc \arrow[bend right=10]{rr}{\scriptstyle \bot}[swap]{\nabla} \&\& \Do \CCc \arrow[bend right=10]{ll}[swap]{\supp}
 \&\&\& \CCc \arrow[bend right=10]{rr}{\scriptstyle \top}[swap]{\Delta} \&\& \Up \CCc \arrow[bend right=10]{ll}[swap]{\inff}
\end{tikzar}
\end{center}
Given a matrix $\Phi:\PObj^\op \times \PSbj\to \Set$, the suprema and the infima \emph{weighted}\/ by its transposes $\Lfun \Phi: \PObj \to \Up \PSbj$ and $\Rfun \Phi: \PSbj \to \Do \PObj$ can similarly be viewed as adjoints:
\begin{center}
\begin{tikzar}{}
\PSbj \arrow[bend right=10]{rr}{\scriptstyle \bot}[swap]{\Rfun \Phi} \&\& \Do \PObj \arrow[bend right=10]{ll}[swap]{\supp_\Phi}
 \&\&\& \PObj \arrow[bend right=10]{rr}{\scriptstyle \top}[swap]{\Lfun \Phi} \&\& \Up \PSbj \arrow[bend right=10]{ll}[swap]{\inff_\Phi}
\end{tikzar}
\end{center}
It is, of course, well known and easy to see that the weighted limits can in ordinary categories be reduced to the ordinary limits. The situation is slightly more subtle with the weighted inferior and superior limits. To align the two situations, note that the adjunctions
\[
\PSbj\left(\sideset{}{_\Phi}\supp  \lft \alpha, b\right) \cong \Do \PObj\left(\lft\alpha, \Rfun\Phi b \right) \qquad \qquad \qquad  \PObj\left(a, \sideset{}{_\Phi}{\inff} \rgt \beta \right) \cong \Up\PSbj\left(\Lfun\Phi a, \rgt \beta \right)
\]
will now become
\[
\PObj\left(\sideset{}{_\Phi}\liminff  \rgt \beta, a\right) \cong \Do \PObj^{\lft\Phi}\left(\ran \Phi \rgt\beta, \nabla_\Phi a \right) \qquad \qquad \qquad  \PSbj\left(b, \sideset{}{_\Phi}{\limsupp} \lft \alpha \right) \cong \Up\PSbj^{\rgt\Phi}\left(\Delta_\Phi b, \lan \Phi \lft \alpha \right)
\]

\begin{definition}\label{def-liminfsup-Phi}
Given a matrix $\Phi: \PObj^\op \times \PSbj \to \Set$, with the induced extensions as in Fig.~\ref{fig-extensions}, we define the operations $\Phi$-\emph{limit inferior} $\liminff_\Phi$ and $\limsupp_\Phi$ by the following adjunctions
\begin{center}
\begin{tikzar}{}
\PObj \arrow[bend right=10]{rr}{\scriptstyle \bot}[swap]{\nabla_\Phi} \&\& \Do \PObj^{\lft\Phi} \arrow[bend right=10]{ll}[swap]{\liminff_\Phi}
 \&\&\& \PSbj \arrow[bend right=10]{rr}{\scriptstyle \top}[swap]{\Delta_\Phi} \&\& \Up \PSbj^{\rgt \Phi} \arrow[bend right=10]{ll}[swap]{\limsupp_\Phi}
\end{tikzar}
\end{center}
where $\nabla_\Phi$ and $\Delta_\Phi$ are
 as defined in Fig.~\ref{fig-extensions}.
\end{definition}

\subsubsection{Two pairs of "Yoneda embeddings"}
In this section we spell out the basic properties of the two kinds of "Yoneda embeddings" induced by a matrix $\Phi:\AAa\biarrow \BBb$:

\begin{itemize}
\item $\lft \Phi$-algebra representables and $\rgt \Phi$-coalgebra representables
\[\nabla_\Phi : \AAa \to \Do \AAa^{\lft \Phi}\qquad\qquad\qquad \Delta_\Phi: \BBb \to \Up \BBb^{\rgt \Phi} \]

\item $\Phi$-representable presheaves and postsheaves
\[\Lfun \Phi : \AAa \to \Up \BBb^{\rgt \Phi}\qquad\qquad\qquad \Rfun \Phi: \BBb \to \Do \AAa^{\lft \Phi} \]
\end{itemize}
The underlying functors are as in Fig.~\ref{fig-extensions}. The structures are as follows.

\begin{lemma} Every presheaf $\lft \alpha \in \Do \PObj$ induces the $\rgt \Phi$-coalgebra $\lan \Phi \lft \alpha \tto{\lan \Phi \eta} \lan \Phi \ran \Phi \lan \Phi \lft \alpha$. Every $\Phi$-representable postsheaf $\Lfun \Phi a$ is thus canonically a $\rgt \Phi$-coalgebra, since $\Lfun \Phi a = \lan \Phi \nabla a$. 

Any postsheaf $\rgt \beta \in \Up \PSbj$ induces the $\lft \Phi$-algebra $\ran \Phi \rgt \beta \oot{\ran \Phi \varepsilon} \ran \Phi \lan \Phi \ran \Phi \rgt \beta$. Every $\Phi$-representable presheaf $\Rfun \Phi b$  is thus canonically a $\lft \Phi$-algebra, since  $\Rfun \Phi b = \ran \Phi \Delta b$.
\end{lemma}


\begin{lemma}[Matrix Yoneda Lemma]\label{lemma-yoneda}
For every $a\in \PObj$ and every $\rgt \beta \in \Up \PSbj^{\rgt \Phi}$ there is  a natural  bijection
\bea\label{eq-Yonedaa}
\Do \PObj^{\lft \Phi} \left(\nabla_\Phi a, \ran \Phi \rgt \beta\right) & \cong & \ran \Phi \rgt \beta (a)
\eea
For every $b\in \PSbj$ and every $\lft \alpha \in \Do \PObj^{\lft \Phi}$ there is  a natural  bijection
\bea\label{eq-Yonedab}
\Up \PSbj^{\rgt \Phi} \left(\lan \Phi \lft \alpha, \Delta_\Phi b \right) & \cong & \lan \Phi \lft \alpha (b)
\eea
\end{lemma}


\begin{corollary}[Matrix Yoneda embedding] $\Do \PObj^{\lft \Phi}\left( \nabla_\Phi a, \Rfun\Phi b\right)\ \cong \ \Phi(a, b) \ \cong \ \Up\PSbj^{\rgt \Phi} \left(\Lfun \Phi a, \Delta_\Phi b \right)$  

\end{corollary}

\subsubsection{Completeness and generation}


\begin{corollary}\label{cor-comple-phi}
$\Do \PObj^{\lft \Phi}$ is $\liminff_\Phi$-complete. $\Up \PSbj^{\rgt \Phi}$ is $\limsupp_\Phi$-complete.  
\end{corollary}

\begin{proposition}\label{prop-represgen-phi}
Every $\lft \Phi$-algebra is a limit inferior in $\Do \AAa^{\lft\Phi}$ of $\lft\Phi$-algebra representables. Every $\rgt\Phi$-coalgebra is a limit superior in $\Up \BBb^{\rgt \Phi}$ of $\rgt\Phi$-coalgebra representables. 
\end{proposition}

\begin{theorem}\label{thm-completion-phi}
The $\Phi$-extended Yoneda embeddings realize the $\liminff_\Phi$-completion and $\limsupp_\Phi$-completion:\begin{itemize}
\item $\nabla_\Phi : \AAa \tto\nabla \Do \AAa\tto{\overline \Phi} \Do \AAa^{\lft \Phi}$ is the $\liminff$-completion of $\CCc$, whereas
\item $\Delta_\Phi : \BBb \tto{\Delta} \Up \BBb\tto{\underline \Phi} \Up \CCc^{\rgt \Phi}$ is the $\limsupp$-completion of $\CCc$.
\end{itemize}
\end{theorem}

\subsection{Minimal bicompletion of a matrix}\label{Sec-two-three}
\subsubsection{Loose extensions}
In general, a matrix $\Phi: \PObj \biarrow \PSbj$ always induces a \emph{loose}\/ extension $\UD \Phi : \Do \PObj^{\lft \Phi} \biarrow \Up \PSbj^{\rgt\Phi}$, defined
\bea\label{eq-loose}
\UD\Phi \left(a, b \right) & = & \left\{
f\in \Do \PObj\left(\lft\alpha, \ran\Phi \rgt \beta\right)\ \mbox{\LARGE $\Bigg|$}\ 
\begin{tikzar}{}
\ran \Phi\lan\Phi \lft \alpha \ar{dd}[swap]{a} \ar{rr}{\ran \Phi\lan\Phi f} \&\& \ran\Phi\lan\Phi\ran\Phi \rgt\beta \\ \\
\lft\alpha \ar{rr}[swap]{f} \&\& \ran\Phi\rgt \beta \ar{uu}[swap]{\ran\Phi b}
\end{tikzar}
\right\}
\eea

\begin{proposition}
Each of the following squares commutes if and only if the other one commutes.  
\[
\begin{tikzar}{}
\ran\Phi\lan\Phi \lft \alpha \ar{dd}[swap]{a} \ar{rr}{\ran\Phi\lan\Phi f} \&\&\ran\Phi\lan\Phi\ran\Phi\rgt\beta \&\& \lan\Phi\ran\Phi\lan\Phi \lft\alpha \ar{dd}[swap]{\lan\Phi a} \ar{rr}{\lan\Phi\ran\Phi f'} \&\&\lan\Phi\ran\Phi\rgt\beta \\ \&\&\& \iff \\
\lft\alpha \ar{rr}[swap]{f} \&\& \ran\Phi \rgt \beta \ar{uu}[swap]{\ran\Phi b} \&\& \lan\Phi\lft \alpha \ar{rr}[swap]{f'} \&\& \beta \ar{uu}[swap]{b}
\end{tikzar}
\]
The commutativity of the preceding squares implies the commutativity of the following squares, which are each other's transposes.
\bear
\begin{tikzar}{}
\&\&\&\& \lan \Phi \ran \Phi \lan \Phi \lft \alpha \ar{r}{\lan \Phi\ran \Phi f'} \& \lan \Phi \ran \Phi \rgt \beta \\
\lft \alpha \arrow{r}{f} \& \ran \Phi \rgt \beta \arrow[yshift = 1ex]{r}{\ran\Phi b} \arrow[yshift = -1ex]{r}[swap]{\eta} \&  \ran\Phi \lan \Phi \ran \Phi  \rgt \beta \& \iff \\
\&\&\&\& \lan \Phi \lft \alpha\ar{uu}{\ran\Phi\eta} \ar{r}[swap]{f'} \& \rgt \beta \ar{uu}[swap]{b}
\end{tikzar} 
\eear
\bear
\begin{tikzar}{}
\&\&\&\& \ran \Phi \lan \Phi \lft \alpha \ar{dd}[swap]{a} \ar{r}{\ran \Phi \lan \Phi f'} \& \ran \Phi  \lan \Phi  \ran \Phi \rgt \beta  \ar{dd}{\ran \Phi \varepsilon}\\
\lan \Phi \ran \Phi \lan \Phi \lft \alpha 
 \arrow[yshift=1ex]{r}{\lan\Phi a} \arrow[yshift=-1ex]{r}[swap]{\varepsilon}\&  
 \lan \Phi \lft \alpha \arrow{r}{f} \& \rgt \beta \& \iff \\
\&\&\&\& \lft \alpha \ar{r}[swap]{f'} \& \ran \Phi \rgt \beta
\end{tikzar} 
\eear
\end{proposition}

\begin{conjecture} $\UD \Phi$ isomorphic with the matrix
\bear \UD\Phi \left(a, b \right) & = & \Do{\left(\PObj\times\PSbj^{\op}\right)}^{\lft\Phi\times \rgt\Phi} \left(\lft \alpha\times \rgt\beta, \Phi \right)
\eear
which is equivalent to the matrix of the adjunction $\Lan\Phi \dashv \Ran\Phi: \PSbj^{\rgt\Phi} \to \Do \PObj^{\lft \Phi}$, defined 
\bea\label{eq-cat-extensionone-phi}
\Lan \Phi \lft \alpha\left(u\right) & = &  \Do \PObj^{\lft\Phi} \left(\lft \alpha, \Lfun{\Phi} u\right)\\  
\Ran \Phi \rgt \beta\left(\ell\right) & = &  \Up \PSbj^{\rgt\Phi}  \left(\Rfun{\Phi} \ell, \rgt \beta\right) \label{eq-cat-extensiontwo-phi}
\eea
with the structure maps induced by composition with the structure maps $a:\ran\Phi\lan\Phi\lft \alpha\to \lft\alpha$ and  $b:\rgt \beta\to \lan\Phi\ran\Phi\rgt\beta$.
\end{conjecture}

\subsubsection{Tight extensions}
But this loose extension is of little semantical value. E.g., when $\Phi$ is a partial ordering like in \eqref{eq-order}, $\UD\Phi$ picks all pairs of a saturated lower set and a saturated upper set which are contained in each other's sets of bounds, but do not necessarily contain \emph{all}\/ such bounds. So it does not capture the Dedekind cuts. 

The tight extension $\cut \Phi$ brings us closer to the Dedekind cuts:
\bea\label{eq-tight}
\cut\Phi \left(a, b \right) & = & \left\{
f\in \UD \Phi\left(a,b\right)\ |\ \mbox{$f$ is mono, and $f'$ is epi}
\right\}
\eea
Since $\Do \PObj^{\lft\Phi}$ and $\Up \PSbj^{\rgt\Phi}$ are regular categories, $\cut \Phi$ can be extracted from $\UD\Phi$ by two closure operators: first extracting the mono factors, and then the epis of their transposes, or equivalently the other way around. After the factorizations, in the first case the transpose of the resulting epi will be mono; in the second the transpose of the resulting mono will be epi. Either way, the process will stop.

The resulting matrix $\cut \Phi$ will be a reflective submatrix of $\UD\Phi$. The completeness and the generation will be inherited, but tight. We need to prove that the inferior limits that existed in $\PObj$ and the superior limits that existed in $\PSbj$ are preserved.

\begin{conjecture}
For every matrix $\Phi: \PObj \to \PSbj$, the tight extension $\cut \Phi: \Do \PObj^{\lft \Phi} \biarrow \Up \PSbj^{\rgt\Phi}$ is the minimal bicompletion.
\end{conjecture}

\section{When does limit inferior boil down to limit?}\label{Sec-three}
By the couniversal property of the (Eilenberg-Moore) categories of algebras for a monad \cite[Part~0.6]{Lambek-Scott:book}, there are always the comparison adjunctions between $\Do\PObj$ and $\Up\PSbj^{\rgt \Phi}$, and between $\Up\PSbj$ and $\Do\PObj^{\lft \Phi}$, as displayed in the leftmost diagram of Fig.~\ref{fig-comp}, since the monad $\lft \Phi$ and the comonad $\rgt \Phi$ are induced by the adjunction $\lan \Phi \dashv \ran \Phi: \Up\PSbj \to \Do \PObj$. When these comparisons are equivalences, then this adjunction transfers to the two Eilenberg-Moore categories, as indicated in the rightmost diagram of Fig.~\ref{fig-comp}. Moreover, the inferior $\Phi$-limits $\limsupp_\Phi$ in $\PSbj$ then boil down to the suprema $\supp$ in $\PObj$, whereas the superior $\Phi$-limits $\limsupp_\Phi$ in $\PObj$ boil down to the infima $\inff$ in $\PSbj$. In terms of the concept mining example from the Introduction, the structural components represented in $\Do\PObj^{\lft \Phi}$ can be computed as infima functions in $\Up\PSbj$, whereas the functional modules represented in $\Up \PSbj^{\rgt \Phi}$ can be computed as suprema of parts in $\Do \PObj$. Connecting the extensions $\UD\Phi$ and $\cut \Phi$ along the equivalences $\Do\PObj\simeq \Up\PSbj^{\rgt \Phi}$ and $\Up\PSbj\simeq \Do\PObj^{\lft \Phi}$ shows that all loose extensions are already tight.
\begin{proposition}
For any matrix $\Phi:\PObj \biarrow \PSbj$, the extensions $\lan \Phi \dashv \ran \Phi: \Up\PSbj \to \Do \PObj$ are both monadic if and only if the loose and the tight extensions coincide, i.e. $\UD \Phi \simeq \cut \Phi$.
\end{proposition}
The notion of \emph{monadicity}\/ \cite[Sec.~VI.7]{CWM} here precisely captures the equivalences of interest, as  $\Do\PObj\simeq \Up\PSbj^{\rgt \Phi}$ means that $\ran \Phi$ is monadic and $\Up\PSbj\simeq \Do\PObj^{\lft \Phi}$ means that $\lan \Phi$ is monadic. In this section, we  study the monadicity of the extensions $\ran \Phi$ and $\lan \Phi$ in order to gain insight into the situations when the inferior and superior limits boil down to the ordinary limits, and the situations when they genuinely provide new information.
\begin{figure}[htbp]
\centering
\begin{tikzar}{row sep=3.6em,column sep=3.6em}
\Do\AAa
\arrow[->,loop above]{}{\Lmon\Phi}
\arrow[->,dashed]{dr}{}
\arrow[->,bend right=15]{d}[swap]{\lan\Phi}
\arrow[<-,bend left=15]{d}{\ran\Phi}
\arrow[no line]{d}{\scriptstyle\dashv}
\arrow[->,bend left=15]{r}{\overline{\Phi}}
\arrow[<-,bend right=15]{r}[swap]{U}
\arrow[no line]{r}{\scriptstyle\bot}
\&
\Do\AAa^{\Lmon\Phi}
\\
\Up\BBb
\arrow[->,loop below]{}{\Rmon\Phi}
\arrow[->,dashed]{ur}{}
\&
\Up\BBb^{\Rmon\Phi}
\arrow[->,bend right=15]{l}[swap]{V}
\arrow[<-,bend left=15]{l}{\underline{\Phi}}
\arrow[no line]{l}{\scriptstyle\bot}
\end{tikzar}
\qquad
\begin{tikzar}{}
\Do\AAa
\arrow[->,loop above]{}{\Lmon\Phi}
\arrow[->,bend right=5,shift right=2pt]{ddd}[swap]{\lan\Phi}
\arrow[<-,bend left=5,shift left=2pt]{ddd}{\ran\Phi}
\arrow[no line]{ddd}{\scriptstyle\dashv}
\arrow[->,bend left=5,shift left=1pt]{rrr}{\overline{\Phi}}
\arrow[<-,bend right=5,shift right=1pt]{rrr}[swap]{U}
\arrow[no line]{rrr}{\scriptstyle\bot}
\&\&\&
\Do\AAa^{\Lmon\Phi}
\\
\&
\CCC
\arrow[->,loop above]{}{\Lmon\Phi|_\CCC}
\arrow[>->]{ul}{}
\arrow[->,dashed]{dr}{}
\arrow[->,bend right=15]{d}[swap]{\lan\Phi|_\CCC}
\arrow[<-,bend left=15]{d}{\ran\Phi|_\DDD}
\arrow[no line]{d}{\scriptstyle\dashv}
\arrow[->,bend left=15,pos=0.6]{r}
\arrow[<-,bend right=15]{r}[swap]{}
\arrow[no line]{r}{\scriptstyle\bot}
\&
\CCC^{\Lmon\Phi|_\CCC}
\arrow[-,shift left=1pt]{ur}
\arrow[-,shift right=1pt]{ur}
\&
\\
\&
\DDD
\arrow[->,loop below]{}{\Rmon\Phi|_\DDD}
\arrow[>->]{dl}{}
\arrow[->,dashed]{ur}{}
\&
\DDD^{\Rmon\Phi|_\DDD}
\arrow[-,shift left=1pt]{dr}
\arrow[-,shift right=1pt]{dr}
\arrow[->,bend right=15]{l}[swap]{}
\arrow[<-,bend left=15,pos=0.4]{l}
\arrow[no line]{l}{\scriptstyle\bot}
\&
\\
\Up\BBb
\arrow[->,loop below]{}{\Rmon\Phi}
\&\&\&
\Up\BBb^{\Rmon\Phi}
\arrow[->,bend right=5,shift right=1pt]{lll}[swap]{V}
\arrow[<-,bend left=5,shift left=1pt]{lll}{\underline{\Phi}}
\arrow[no line]{lll}{\scriptstyle\bot}
\end{tikzar}
\qquad
\begin{tikzar}{row sep=3.6em,column sep=3.6em}
\Do\AAa
\arrow[->,loop above]{}{\Lmon\Phi}
\arrow[->,bend right=15]{d}[swap]{\lan\Phi}
\arrow[<-,bend left=15]{d}{\ran\Phi}
\arrow[no line]{d}{\scriptstyle\dashv}
\&
\Do\AAa^{\Lmon\Phi}
\arrow[>->]{dl}{}
\\
\Up\BBb
\arrow[->,loop below]{}{\Rmon\Phi}
\&
\Up\BBb^{\Rmon\Phi}
\arrow[>->]{ul}{}
\arrow[->,bend left=15]{u}{}
\arrow[<-,bend right=15]{u}[swap]{}
\arrow[no line]{u}{\scriptstyle\dashv}
\end{tikzar}

\caption{Comparisons between the $\protect\inff$- and $\protect\liminff_\Phi$-completions, and between the $\protect\supp$- and $\protect\limsupp_\Phi$-completions}
\label{fig-comp}
\end{figure}

\subsection{Monadicity workflow}\label{Sec-three-one}
%
As a reminder, we quote the Precise Monadicity Theorem in Appendix B. Intuitively, its impact on the concrete instances of our situation is that it allows constructing the inferior limits, which are in principle the suprema of infima, as specific maximal cones into the infima.

We begin describing a convenient setting of subcategories, as displayed in the middle in Fig.~\ref{fig-comp}. When  $\Phi_*\colon\Up\BBb\to\Do\AAa$ restricts to
a monadic functor $\DDD\to\CCC$, so that $\DDD\simeq \CCC^{\Lmon\Phi|_\CCC}$,
then  we have an embedding $(\Do\AAa)^{\Lmon\Phi}\mono\Up\BBb$. 

In the general framework of an adjunction as in Appendix B, items (a-b) of the Monadicity Theorem say that the induced Eilenberg-Moore category $\CCC^T$ is coreflective within the category $\DDD$ whenever $\DDD$ has and $U$ preserves reflexive $U$-split coequalizers.
However, its converse does not hold%
. The task is thus to spell out the full subcategories
$\CCC\subseteq\Do\AAa,\, \DDD\subseteq\Up\BBb$
explicitly. 
Towards this goal, and to simplify calculations with the algebras, we propose the following.
%

\begin{definition}
  An object $B$ is said to be a retract of an object $A$
  if there exist morphisms $B\to A\to B$ whose composite
  is $\id_B$.
  For a full subcategory $\AAA\subseteq\EEE$,
  we denote by $\Retr_\EEE(\AAA)\subseteq\EEE$
  the full subcategory of all retracts in $\EEE$
  of objects in $\AAA$.
\end{definition}

\paragraph{Notational conventions.}
  For a functor $G$, we denote its full image by $\Img G$.
  For a category $\EEE$ and its full subcategories $\AAA,\AAA'$,
  we loosely use $\AAA\subseteq\AAA'$ to denote
  any object in $\AAA$ is isomorphic in $\EEE$ to some object in $\AAA'$.

\begin{lemma}\label{lemma:retrAdj}
  Let $F\dashv U\colon\YYY\to\XXX$ be an adjunction.
  \begin{enumerate}
    \item
      Let $\CCC\subseteq\XXX$ be a full subcategory.
      If $\Retr_\XXX(\Img U)\subseteq\CCC$,
      there exists a canonical equivalence of categories
      $\XXX^T \simeq \CCC^{T|_{\CCC}}$.
    \item
      Assume that $\YYY$ has reflexive $U$-split coequalizers and that
      the left adjoint $L\colon \XXX^T\to\YYY$ of the comparison functor
      is full and faithful.
      (In particular, 
      $\XXX^T$ is equivalent to a coreflective subcategory $\Img L\subseteq \YYY$.)
      Then,
      $\Retr_\YYY(\Img F)\subseteq \DDD$.
  \end{enumerate}
\end{lemma}

\begin{figure}[htbp]
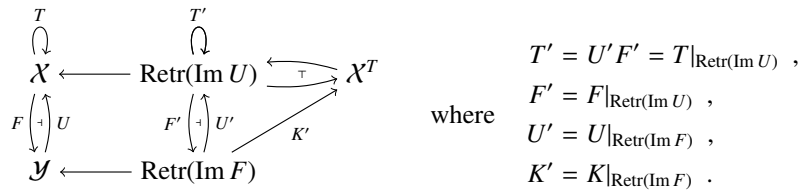

  \centering
  \begin{tikzar}{}
  \XXX
\arrow[->,loop above]{}{T}
\arrow[->,bend right=15]{d}[swap]{F}
\arrow[<-,bend left=15]{d}{U}
\arrow[no line]{d}{\scriptstyle\dashv}
\arrow[<-]{r}{}
  \&
  \Retr(\Img U)
\arrow[->,loop above]{}{T'}
\arrow[->,bend right=15]{d}[swap]{F'}
\arrow[<-,bend left=15]{d}{U'}
\arrow[no line]{d}{\scriptstyle\dashv}
\arrow[->,loop above]{}{}
\arrow[->,bend right=10]{r}[swap]{}
\arrow[<-,bend left=10]{r}{}
\arrow[no line]{r}{\scriptstyle\top}
  \&
  \XXX^T
  \\
  \YYY
\arrow[<-]{r}{}
  \&
  \Retr(\Img F)
\arrow[->]{ur}[swap]{K'}
  \end{tikzar}
  \quad where\quad
  $\begin{aligned}
  T'&=U'F'=T|_{\Retr(\Img U)}
  \enspace,\\
  F'&=F|_{\Retr(\Img U)}
  \enspace,\\
  U'&=U|_{\Retr(\Img F)}
  \enspace,\\
  K'&=K|_{\Retr(\Img F)}
  \enspace.
  \end{aligned}$
  \caption{Restricting an adjunction to retracts of images}
  \label{fig:retim}
\end{figure}

The above lemma intuitively means
\begin{itemize}
  \item we need at most retracts of images under $U$ in $\XXX$, and that
  \item we need at least retracts of images under $F$ in $\YYY$,
\end{itemize}
in order to obtain a monadic functor of the form
$\XXX^T\simeq\DDD\tto{U|_{\DDD}}\CCC$.
In the later discussion, we restrict an adjunction as Fig.~\ref{fig:retim}
and calculate the category of $T$-algebras by $\XXX^T = (\Retr(\Img U))^{T'}$.

\subsection{Completing constant matrices}\label{Sec-three-two}
Any set $R$ can be viewed as a constant matrix $\widetilde{R}:1\biarrow 1$ by setting $\widetilde{R}(0) = R$, where $1=\{0\}$. We abuse notation and write $\widetilde{R}$ as $R$. The extensions $\lan R \dashv \ran R \colon\Set^\op\to\Set$ are thus $\lan R X = \ran R X = R^X$, and they induce the monad $\lft R X = R^{R^X}$ on $\Set$, and the same comonad $\rgt R$ on $\Set^\op$.

Lemma~\ref{lemma:weakReflEq} in the Appendix B  helps characterizing the monadicity of $\lan R$ and $\ran R$.

\begin{proposition} For a set $R$ with at least 2 elements, the functor $\widetilde{R}_\ast \colon\Set^\op\to\Set$  is monadic. When $R$ is a singleton, then the monad $\lft  R:\Set \to \Set$ has a single algebra, and the comonad $\rgt R:\Set^\op \to \Set^\op$ has a single coalgebra.  When $R$ is empty, then the they have two algebras and coalgebras respectively.
\end{proposition}

\begin{corollary}
The loose extension of the constant matrix $R$ is always in the form $\UD R: \Set\biarrow  \Set$ with $\UD R(X,Y) = \Set(X,Y)$. The tight extension is
\begin{itemize}
\item $\cut R = \UD R:\Set \biarrow \Set$ when $R$ has at least 2 elements
\item $\cut 1 : 1 \biarrow 1$ with $\cut 1(0, 0) = 1$, where $1 =\{0\}$
\item $\cut 0 : 2 \biarrow 2$ with $ \cut 0 (x,y) = 1$ if and only if $x\leq y$ within $2 = \{0\lt 1\}$.
\end{itemize}
\end{corollary}

\subsection{Completing groups}\label{Sec-three-three}
Let $\CCc$ be a group $G$, viewed as a one object category with invertible morphisms.
The category $\Do G$ of presheaves is the category of right $G$-sets,
or the category $\Seta{G^\op}$ of (left) $G^\op$-sets.
Indeed as a discrete fibration over $1$, the total category of the presheaf is
a set $X$ with an action $X\times G\to X$. The adjunction $H^* \dashv H_*$ is given explicitly as follows.
We think of $G$ as a (left $G$, right $G$)-set by the multiplication.
For a right $G$-set $X$, the (left) $G$-set $H^* X$ is
the set $\Seta{G^\op}(X,G)$ with the action $(g\cdot f)(x) = g(f(x))$.
Similarly, $H_* Y=\Seta{G}(Y,G)$ for a left $G$-set $Y$.

\begin{proposition}\label{proposition:retrImGSet}
  We have
  $\Img H^* \simeq \{0\} \cup \{G^I \mid I\in\Set\}$
  and
  $\Retr(\Img H^*) \simeq \{1\} \cup \{G\times I \mid I\in\Set\}$,
  where $G^I$ is the exponential in $\Set$ with the pointwise multiplication
  $(g\cdot f)(i) = g(f(i))$
  and
  $G\times I$ is the free $G$-set generated by the set $I$
  (i.e.\ $g\cdot(h,i)=(gh,i)$).
\end{proposition}

We denote by $\Setaf{G}$ the full subcategory
$\{1\} \cup \{G\times I \mid I\in\Set\} \subseteq \Seta{G}$
of a singleton and free $G$-sets.

\begin{proposition}\label{proposition:monadicGSet}
  The functor $H_*\colon(\Setaf{G})^\op\to \Setaf{G^\op}$ is monadic. 
\end{proposition}

\begin{proposition}
The category $(\Do G)^{\LRC}$ of algebras are equivalent with the full subcategory of $\Up G$ which consists of free left $G$-sets, together with the singleton.
\end{proposition}

\begin{corollary}
The loose extension of a nontrivial group is the canonical connection of its left and right actions. The tight extension is the canonical extension of its free actions.
\end{corollary}

\subsection{Completing posets}\label{Sec-three-four}
Let $\CCc$ be a poset $(P,\le)$. We write the poset of lower sets of $P$ as $\Lower{P}$, and the poset of upper sets of $P$ by $\Upper{P}$. They are respectively the join and the meet completions. While $P$'s categorical supremum completion $\Do P = \Set^{P^\op}$ and its  infimum completion $\Up P =\left(\Set ^P\right)^\op$ are proper categories, its limit inferior completion $\Do P^{\LRC}$, and its limit superior completion $\Up P^{\RLC}$, although still constructed over $\Set$ --- turn out to be both equivalent to a lattice, and in particular to $P$'s Dedekind-MacNeille completion $\UpperLower P$.

\begin{lemma} The lattice of subobjects of the terminal object in $\Do P$ is isomorphic to $\Lower P$. The lattice of subobjects of the terminal object of $\Up P$ is isomorphic to $\Upper P$.  In particular, $\Lower P\subseteq \Do P$ and $\Upper P\subseteq \Up P$ are
  full subcategories containing the representables.
\end{lemma}

\begin{lemma}\label{lemma:posetAdjIncl}
  The adjunction $H^*\dashv H_*\colon \Up P \to \Do P$ restricts to
  an adjunction (a Galois connection) between posets $\Lower P, \Upper P$,
  which coincides with $\{0,1\}$-enriched construction.
  Moreover,
  $\Img H^* = \Img (H^*|_{\Lower P})$ and
  $\Img H_* = \Img (H_*|_{\Upper P})$.
\end{lemma}

\begin{corollary}\label{corollary:retrImPoset}
  It holds
  $\Retr_{\Up P}(\Img H^*)=\Img (H^*|_{\Lower P})$.
\end{corollary}

Therefore, the category $(\Do P)^{\LRC}$
is nothing more than
the category of algebras for the adjunction
$\Upper P \leftrightarrows \Lower P$.

\begin{corollary}
  There exist equivalences of categories
  $(\Do P)^\LRC \simeq \UpperLower P \simeq (\Up P)^\RLC$.
\end{corollary}

\begin{corollary}
The tight extension $\cut P$ of a poset $P$ coincides with its Dedekind-MacNeille completion $\UpperLower P$.
\end{corollary}

\subsection{Completing a $\zpz$-vector}\label{Sec-three-five}
A \emph{vector}\/ is a matrix in the form $\Phi : 1 \biarrow \PSbj$. We consider the vectors in $\PSbj = \zpz$, viewed as an additive cyclic group of prime order $p$. Every left $\zpz$-set $X$ has an orbit-decomposition $X = 1\times X_1 + \zpz\times X_p$ where
the action on $1$ is trivial and $\zpz$ is by left multiplication.
We abbreviate this decomposition as
$X = \od{X_1}{X_p}$.

\begin{lemma}
  $ \Seta{\zpz}(\od{X_1}{X_p},\od{Y_1}{Y_p}) \cong
  Y_1^{X_1} (Y_1+pY_p)^{X_p} $.
\end{lemma}
In particular, the adjunction $\Phi^*\dashv\Phi_*$ is explicitly
$\Phi^*(L)=(\od{\Phi_1}{\Phi_p})^L$ and
$\Phi_*(\od{U_1}{U_p})
\cong \Phi_1^{U_1} (\Phi_1+p\Phi_p)^{U_p}$.

\begin{lemma}\label{lemma:reflEqZpZ}
  A reflexive pair $f,g\colon U\rightrightarrows U'$ in $\Seta{\zpz}$
  has some isomorphisms to the bottom row of the diagram
  \begin{center}
    \begin{tikzar}{}
      U \arrow[->]{d}{\cong}
      \arrow[->,shift left=3pt]{r}{f}
      \arrow[->,shift right=3pt]{r}[swap]{g}
      \&
      U' \arrow[->]{d}{\cong}
      \\
      \od{U_1}{U_p}
      \arrow[->,shift left=3pt]{r}{\od{f_1}{f_p}}
      \arrow[->,shift right=3pt]{r}[swap]{\od{g_1}{g_p}}
      \&
      \od{U'_1}{U'_p}
    \end{tikzar}
    \quad for some maps 
    $f_1,g_1\colon U_1\rightrightarrows U'_1,\;
    f_p,g_p\colon U_p\rightrightarrows U'_p$.
  \end{center}

  Let $E_1\tto{e_1} U_1$, $E_p\tto{e_p} U_p$ be equalizers of
  $(f_1,g_1)$, $(f_p,g_p)$ respectively.
  These coequalizer satisfies the condition of Lem.~\ref{lemma:weakReflEq}.2,
  i.e.\ the following diagrams are pullbacks of injections.
\begin{center}\begin{tikzar}{}
  E_1 \arrow[>->]{d}[swap]{e_1} \arrow[>->]{r}{e_1}\& U_1 \arrow[>->]{d}{f_1}
  \&        E_p \arrow[>->]{d}[swap]{e_p} \arrow[>->]{r}{e_p}\& U_p \arrow[>->]{d}{f_p}
  \\
  U_1 \arrow[>->]{r}[swap]{g_1}\& U'_1
  \&        U_p \arrow[>->]{r}[swap]{g_p}\& U'_p
\end{tikzar}\end{center}
\end{lemma}

Let us find
$\Retr(\Img\Phi_*)\subseteq\CCC\subseteq\Set$ and
$\Retr(\Img\Phi^*)\subseteq\DDD\subseteq(\Seta{\zpz})^\op$
to fit the scheme of Fig.~\ref{fig-comp}.
%
\begin{proposition}\label{proposition:monadicZpZ}
  In the following restriction of
  $\Phi^*\dashv\Phi_*\colon(\Seta{\zpz})^\op\to\Set$ makes
  both $(\Phi^*)^\op$ and $\Phi_*$ monadic,
  without changing the categories of algebras:
  \[ \Set^{\Lmon\Phi} = \CCC^{\Lmon\Phi|_{\CCC}} \simeq \DDD
  \enspace,\qquad
  (\Seta{\zpz}^\op)^{\Rmon\Phi} = \DDD^{\Rmon\Phi|_{\DDD}} \simeq \CCC
  \enspace.\]
  \begin{center}
  \begin{tabular}{|c|c|c|}
    \hline
    $\Phi = \od{\Phi_1}{\Phi_p}$
    & $\Phi_p=0$ & $\Phi_p\ge 1$
    \\ \hline
    $\Phi_1=0$
    &
    \begin{tikzar}{}
    \{0,1\}
    \arrow[->,bend right=15]{d}[swap]{} \arrow[<-,bend left=15]{d}{} \arrow[no line]{d}{\scriptstyle\dashv}
    \\
    \{\odx,\odl{}\}^\op
    \end{tikzar}
    &
    \begin{tikzar}{}
    \Set
    \arrow[->,bend right=15]{d}[swap]{} \arrow[<-,bend left=15]{d}{} \arrow[no line]{d}{\scriptstyle\dashv}
    \\
    (\{\odl{}\}\cup\{\odr{U_p}\mid U_p\in\Set\})^\op
    \end{tikzar}
    \\ \hline
    $\Phi_1=1$
    &
    \begin{tikzar}{}
    \{1\}
    \arrow[->,bend right=15]{d}[swap]{} \arrow[<-,bend left=15]{d}{} \arrow[no line]{d}{\scriptstyle\dashv}
    \\
    \{\odl{}\}^\op
    \end{tikzar}
    &
    \begin{tikzar}{}
    \Set
    \arrow[->,bend right=15]{d}[swap]{} \arrow[<-,bend left=15]{d}{} \arrow[no line]{d}{\scriptstyle\dashv}
    \\
    \{\od{}{U_p}\mid U_p\in\Set\}^\op
    \end{tikzar}
    \\ \hline
    $\Phi_1\ge 2$
    &
    \begin{tikzar}{}
    \Set
    \arrow[->,bend right=15]{d}[swap]{} \arrow[<-,bend left=15]{d}{} \arrow[no line]{d}{\scriptstyle\dashv}
    \\
    \{\odl{U_1}\mid U_1\in\Set\}^\op
    \end{tikzar}
    &
    \begin{tikzar}{}
    \Set
    \arrow[->,bend right=15]{d}[swap]{} \arrow[<-,bend left=15]{d}{} \arrow[no line]{d}{\scriptstyle\dashv}
    \\
    (\Seta{\zpz})^\op
    \end{tikzar}
    \\ \hline
  \end{tabular}
  \end{center}
\end{proposition}

\section{Conclusion}\label{Sec-four}

\begin{figure}[htbp]
\begin{center}
\includegraphics[height=3.8cm]{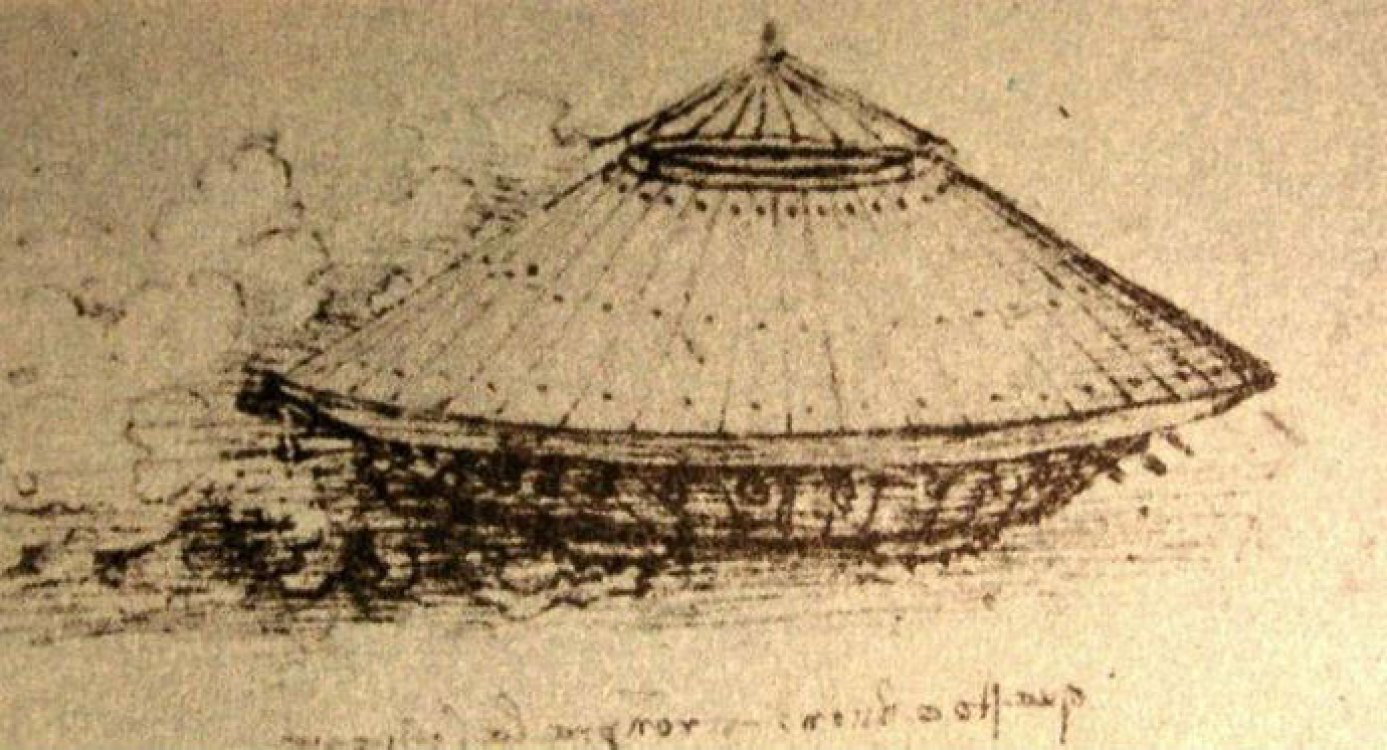}\ \ \ \includegraphics[height=3.8cm]{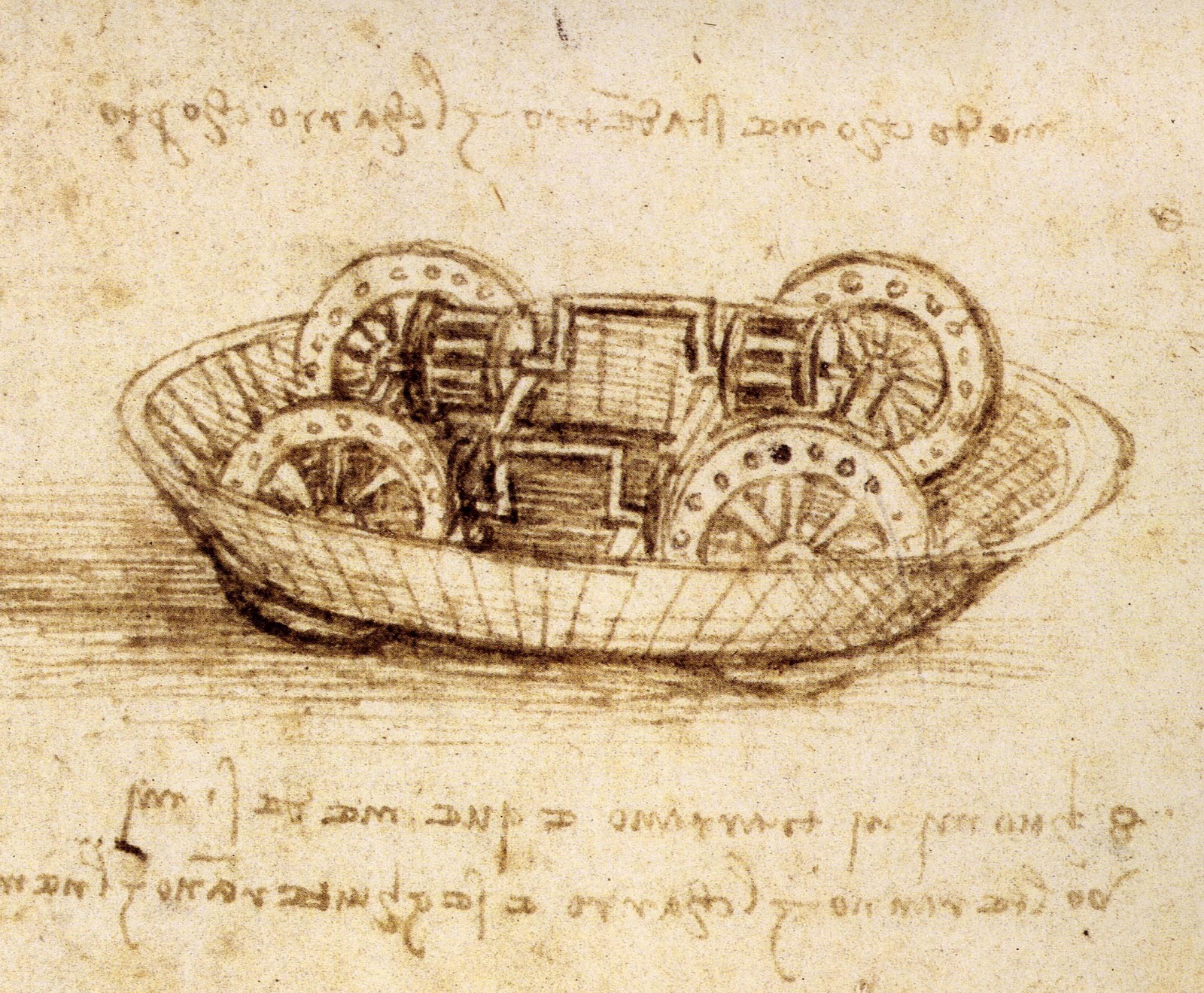}
\caption{Identified object: The external and the internal view}
\label{figtwo}
\end{center}
\end{figure}

Deploying the categorical concept analysis of the unidentified object from Fig.~\ref{figone} according to the technical recipes proposed in this paper, our diligent reader has surely uncovered that the mysterious device consists of two main structural components: the internal mechanism of wheels and gears, and  the external protection shell. On the other hand, the detailed categorical analysis has surely displayed three main functional modules: moving, defending from the outside attacks, and attacking from inside. As desired, the tight matrix then clearly shows that the object must be a model of a man-powered armored combat vehicle from XV century. It was conceived by Leonardo da Vinci, whose drawings are reproduced on Fig.~\ref{figtwo}. The advances of category theory will undoubtedly permit us to better understand Leonardo's conceptualizations of warfare. 

%


\bibliographystyle{plain}
\bibliography{infsup-ref,CT,PavlovicD}

\appendix
\section{Appendix: Proofs}
\begin{prf}{ of  Thm.~\ref{thm-completion}}
Suppose that $\XXx$ is a category with all limits inferior, and that $G:\CCc \to \XXx$ is an arbitrary functor. We show that $G$ has a unique extension $G':\Do \CCc^{\LRC } \to \XXx$, such that
\bea\label{eq-req}
G & = & G'\circ \overline H \circ \nabla
\eea
where $\overline H: \Do \CCc \to \Do \CCc^{\LRC }$, as defined on Fig.~\ref{fig-extensions} instantiated to $\Phi = H$, maps $\CCc$-presheaves to free $\LRC $-algebras, i.e. it is the left adjoint of the forgetful functor $U:\Do \CCc^{\LRC } \to \Do \CCc$.
The construction is illustrated on the following diagram. 
\begin{center}
\begin{tikzar}{row sep = tiny}
\& \& \Do \CCc^{\LRC } \arrow[dashed]{dddd}{G'} \&\& \lft \gamma \oot{\ \ a\ \ } \LRC  \lft \gamma \arrow[mapsto]{dddd}{G'}
\\
\& \Do \CCc \arrow{ru}{\overline H} 
\\
\CCc \arrow{ru}{\nabla} \arrow{rrdd}[swap]{G} 
\\
\\
\& \& \XXx \&\& \liminff\ G\circ\rgt \varphi
\end{tikzar}
\end{center}
Given an arbitrary $\LRC $-algebra $\lft \gamma \oot{\ \ a\ \ } \LRC  \lft \gamma$ in $\Do \CCc$, we construct the equalizer of postsheaves
\begin{center}
\begin{tikzar}{}
\rgt \varphi \arrow[tail]{r} \& \RC  \lft \gamma \arrow[shift left = 1ex]{r}{\RC a} \arrow[shift right = 1ex]{r}[swap]{\eta_{\RC  \lft \gamma}} \&  \RC \LRC  \lft \gamma
\end{tikzar}
\end{center}
which is a coequalizer in $\Up \CCc$. Note that the $\LRC$-algebra $a$ displays the presheaf $\lft \gamma$ as the coequalizer of the $\LC$-image of the pair $\left<\RC a, \eta_{\RC\lft \gamma}\right>$. The $\LRC$-algebra $a$ itself is the coequalizer of the free $\LRC$-algebras over this image, \emph{and}\/ as the limit inferior as decomposed in Proposition~\ref{prop-represgen}. We set the $G'$-image of the $\LRC$-algebra $a$ to be the limit inferior of the functor $\FFf \tto \varphi \CCc \tto G \XXx$. Equation \eqref{eq-req} follows from Lemma~\ref{lemma-repres}. The fact that $G'$ preserves inferior limits follows from the fact that every inferior limit cone $\LC \rgt F \tto{\ \lambda\ } \lft \gamma$
factors through any structure map $\LRC \lft \gamma \tto{ \ a\ } \lft \gamma$: the factorization is the composite $\LC \rgt F \tto{\ \lambda\ } \lft \gamma\tto {\ \eta \ } \LC \lft \gamma$, which obviously boils down to $\lambda$ when further postcomposed with $a$. The uniqueness follows from Proposition~\ref{prop-represgen}.
\end{prf}

\begin{prf}{ of Lemma~\ref{lemma-yoneda}}
Consider a natural transformation $\psi \in \Do \PObj^{\lft \Phi} \left(\ran \Phi \Lfun \Phi a, \ran \Phi \rgt \beta\right)$ on the left-hand side of \eqref{eq-Yonedaa}. By \eqref{eq-cat-extensiontwo} 
and by the naturality of $\psi$, for $f\in \PObj(x,a)$ the left-hand square diagram in Fig.~\ref{figthree} must commute.
\begin{figure}[htbp]
\begin{center}
\begin{tikzar}{column sep = tiny}
\ran \Phi \Lfun \Phi a (a) 
\arrow[-,shift left=1pt]{r} \arrow[-,shift right=1pt]{r}
\ar{dd}{\ran \Phi \Lfun \Phi a(f)} \& 
\Up \PSbj \left( \Lfun \Phi a, \Lfun \Phi a\right) \arrow{rr}{\psi_a} \ar{dd}[swap]{(-)\circ f}\& \& \Up \PSbj \left(\Lfun \Phi a, \rgt \beta \right) 
\arrow[-,shift left=1pt]{r} \arrow[-,shift right=1pt]{r}
\ar{dd}{(-)\circ f}\& 
\ran \Phi \rgt \beta (a)  \ar{dd}[swap]{\ran \Phi \rgt \beta (f)} \&\& \clsr{\nabla_\Phi a}\left(\nabla_\Phi a\right) \arrow{rr}{[\psi]
} \ar{dd}{\clsr{\nabla_\Phi a}\widehat f}\&\& \clsr{\ran \Phi \rgt \beta}\left(\nabla_\Phi a\right)
\ar{dd}[swap]{\clsr{\ran\Phi \rgt \beta}\widehat f}
\\
\\
\ran \Phi \Lfun \Phi a (x) 
\arrow[-,shift left=1pt]{r} \arrow[-,shift right=1pt]{r}
\& 
\Up \PSbj \left( \Lfun \Phi x, \Lfun \Phi a\right) \arrow{rr}{\psi_x} \& \& \Up \PSbj \left(\Lfun \Phi x, \rgt \beta \right) 
\arrow[-,shift left=1pt]{r} \arrow[-,shift right=1pt]{r}
\& 
\ran \Phi \rgt \beta (x)  \&\&\clsr{\nabla_\Phi a}\left(\nabla_\Phi x\right) \arrow{rr}[swap]{[\psi]
} \&\& \clsr{\ran \Phi \rgt \beta}\left(\nabla_\Phi x\right)
\end{tikzar} 
\caption{Matrix Yoneda squares}
\label{figthree}
\end{center}
\end{figure}
Recall that $\lft \Phi$-algebras like $\nabla_\Phi a, \ran \Phi \rgt \beta : \PObj^\op\to \Set$ always canonically extend to functors
$\clsr{\nabla_\Phi a}, \clsr{\ran \Phi \rgt \beta}: \Do \PObj_{\lft \Phi}^\op\to \Set$,
and that $\lft \Phi$-algebra homomorphism $\psi: \nabla_\Phi a\to  \ran \Phi \rgt \beta$ extend to $\clsr{\psi}: \clsr{\nabla_\Phi a} \to \clsr{\ran \Phi \rgt \beta}$. It follows that a $\lft \Phi$-algebra homomorphism $\psi$ must be natural with respect to all homomorphisms between free $\lft \Phi$-algebras, and not just with respect to those arising from $\PObj$. In particular, consider the natural isomorphism
\beq \label{eq-iso}
\Up \PSbj \left( \Lfun \Phi x, \Lfun \Phi a\right)\ \ \stackrel{(a)}=\ \  \ran \Phi \Lfun \Phi a(x)\ \  \stackrel{(b)}\cong\ \  \Do \PObj\left(\nabla x , \ran \Phi \Lfun\Phi a\right)\ \  \stackrel{(c)}=
\ \  \Do \PObj^\op_{\lft \Phi} \left( \nabla_\Phi a,\nabla_\Phi x\right)\eeq
where (a) is based on \eqref{eq-cat-extensiontwo}, (b) on the usual Yoneda lemma, and (c) on the definition of the Kleisli category $\Do \PObj_{\lft \Phi}$. Every $f\in \ \Up \PSbj \left( \Lfun \Phi x, \Lfun \Phi a\right)$ thus induces a unique homomorphism $\widehat f\in \Do \PObj^\op_{\lft \Phi} \left( \nabla_\Phi a,\nabla_\Phi x\right)$, and vice versa.  The naturality condition on $\clsr{\psi}$ now implies that the right-hand square on Fig.~\ref{figthree} must commute, which implies
\beq
\clsr{\psi}_{\nabla_\Phi x}\left( \widehat f \right) \ \ =\ \ \clsr{\psi}_{\nabla_\Phi x}\circ \clsr{\nabla a_\Phi} \widehat f\left(\id_{ \nabla_\Phi a}\right)\ \ = \ \ \clsr{\ran \Phi \rgt \beta}\widehat f\circ \clsr{\psi}_{\nabla_\Phi a} \left(\id_{ \nabla_\Phi a} \right)\ \ =\ \  \clsr{\ran \Phi \rgt \beta}\widehat f \left(\clsr \Psi \right)
\eeq
where $\clsr \Psi = \clsr{\psi}_{\nabla_\Phi a}\left(\id_{ \nabla_\Phi a} \right)$. Hence the bijection between the natural transformations $\clsr{\psi}: \clsr{\nabla_\Phi a} \to \clsr{\ran \Phi \rgt \beta}$ and the elements $\clsr \Psi$ of $\clsr{\ran \Phi \rgt \beta}\left(\nabla_\Phi a\right)$. The restriction to $\psi: \nabla_\Phi a \to \ran \Phi \rgt \beta$ of $\clsr {\psi}$ must be coherent with respect to the natural bijection \eqref{eq-iso}, which means that $\psi$ must be natural with respect to $f \in \Up \PSbj \left( \Lfun \Phi x, \Lfun \Phi a\right)$ just like $\clsr \psi$ was with respect to $\widehat f \in \Do \PObj^\op_{\lft \Phi} \left( \nabla_\Phi a,\nabla_\Phi x\right)$. The naturality of the left-hand square in Fig.~\ref{figthree} now gives
\beq
\psi_x\left(f \right) \ \ =\ \ \psi_x\circ \nabla a_\Phi  f\left(\id_{\Lfun \Phi a}\right)\ \ = \ \ \ran \Phi \rgt \beta( f)\circ \psi_a \left(\id_{\Lfun \Phi a} \right)\ \ =\ \  \ran \Phi \rgt \beta(f) \Psi
\eeq
where $\Psi = \psi_{a}\left(\id_{\Lfun \Phi a} \right)$. Hence the bijection between the $\lft \Phi$-algebra homomorphisms $\psi\in \Do\PObj^{\lft \Phi}\left( \nabla_\Phi a, \ran \Phi \rgt \beta\right)$ and the elements $\Psi$ of $\ran \Phi \rgt \beta\left(a\right)$, as claimed in \eqref{eq-Yonedaa}. Claim \eqref{eq-Yonedab} is proven dually.
\end{prf}

\begin{prf}{ of Lem.~\ref{lemma:retrAdj}}
  \begin{enumerate}
    \item
      Let $a\colon UFC\to C$ be a $T$-algebra.
      By the unit law of Eilenberg-Moore algebras,
      we have a retract $C\overset{\eta_C}\mono UFC\overset{a}\epi C$.
      In particular, the underlying object $C$ of the algebra is
      a retract of an image under $U$.
    \item
      Let $\DDD'=\Img L$.
      Let $A$ be an object of $\DDD'$, and
      $B\mono A\epi B$ be a retract of $A$ in $\DDD$.
      By the monadicity, the upper row of the following diagram is a coequalizer
      in $\DDD'$.
      \begin{center}
      \begin{tikzar}{}
      A
      \&\arrow[->]{l}{\varepsilon A}
      FUA
      \&\arrow[->,shift right=3pt]{l}[swap]{\varepsilon FUA}
      \arrow[->,shift left=3pt]{l}{FU\varepsilon A}
      FUFUA
      \\
      B
      \arrow[->,shift left=3pt]{u}{s}
      \arrow[<-,shift right=3pt]{u}[swap]{r}
      \&\arrow[->]{l}{\varepsilon B}
      FUB
      \arrow[->,shift left=3pt]{u}{FUs}
      \arrow[<-,shift right=3pt]{u}[swap]{FUr}
      \&\arrow[->,shift right=3pt]{l}[swap]{\varepsilon FUB}
      \arrow[->,shift left=3pt]{l}{FU\varepsilon B}
      FUFUB
      \arrow[->,shift left=3pt]{u}{FUFUs}
      \arrow[<-,shift right=3pt]{u}[swap]{FUFUr}
      \end{tikzar}{}
      \quad where $r\circ s=\id_B$
      \end{center}
      Hence the upper row is a coequalizer in $\DDD$ since
      $\DDD'\subseteq\DDD$ is a coreflective subcategory.
      The squares commutes serially, and
      all the columns are retracts.
      It is a straightforward consequence that
      the lower row of the diagram is a coequalizer.
      Equivalently, the counit $LKB\to B$ of $L\dashv K$ at $B$
      is an isomorphism, and thus $B\in\DDD'$.
      \vspace{-\baselineskip}
  \end{enumerate}
\end{prf}


\begin{prf}{ of Prop.~\ref{proposition:retrImGSet}}
  Let $X$ be a right $G$-set.
  If there exists a right $G$-map,
  there exists an isomorphism $X\cong I\times G$ for some set $I$
  by Lem.~\ref{lemma:repnFree}.
  Then, we have $H^*X = \Seta{G^\op}(X,G) \cong \Seta{G^\op}(I\times G, G)
  \cong \Set(I,G) = G^I$.  If there exist no right $G$-maps, for instance $X=1$,
  we have $H^*X = 0$.
%
\end{prf}

\begin{lemma}\label{lemma:repnFree}
  Let $X$ be a $G$-set and $J$ be a set.
  A $G$-map $f\colon X\to G\times J$ to the free $G$-set generated by $J$
  is a composite
  $X\cong G\times I\tto{\id_G\times k} G\times J$
  for some $k\colon I\to J$ in $\Set$.
\end{lemma}

\begin{prf}{ of Lemma~\ref{lemma:repnFree}}
  Let $I=f^{-1}(\{e\}\times J)$.
  The action of $X$ induces an isomorphism $X\cong G\times I$.
\end{prf}

\begin{lemma}\label{lemma:retrGSet}
  A retract of a singleton in $\Seta{G}$ is a singleton.
  A retract of a free $G$-set is free.
\end{lemma}
\begin{prf}{ of  Lemma~\ref{lemma:retrGSet}}
  The first claim is obvious.
  The latter claim is by Lem.~\ref{lemma:repnFree}.
\end{prf}

\begin{prf}{ of  Prop.~\ref{proposition:monadicGSet}}
  By the Monadicity Theorem, 
  this proposition reduces to the following two lemmas.
\end{prf}
\begin{lemma}
  The following holds.
  \begin{enumerate}
    \item The category $\Setaf{G}$ has reflexive equalizers.
    \item The functor $H_*\colon(\Setaf{G})^\op\to \Setaf{G^\op}$
      preserves reflexive coequalizers.
  \end{enumerate}
\end{lemma}
\begin{prf}{} 
  By Lem.~\ref{lemma:retrGSet}, a reflexive pair in $\Setaf{G}$ is either
  $1\rightrightarrows 1$ or $G\times I\rightrightarrows G\times J$.
  The pair $1\rightrightarrows 1$ trivially has an equalizer
  that is preserved by any functor.

  Let $r\colon G\times J\to G\times I$ be a common retraction in $\Setaf{G}$
  of the pair $(f,h)\colon G\times I\rightrightarrows G\times J$.
  We may assume $r=\id_G\times r'$ for some map $r'\colon J\to I$ by
  Lem.~\ref{lemma:repnFree}.
  Define a map $f'\colon I\to J$ by
  $\bigl(g,f'(i)\bigr)=f(e,i)$ for each $i\in I$ where
  we have $g=e$ by $(e,i)=r(f(e,i))=(g,r'(f'(i)))$.
  Moreover for any $g'\in G$,
  $f(g',i) = f(g'(e,i)) = g'(f(e,i)) = g'(e,f'(i)) = (g',f'(i))$
  Therefore, there exists maps $(f',h')\colon I\rightrightarrows J$ such that
  $f=\id_G\times f'$, $h=\id_G\times h'$, and
  $r'$ is a common retraction of $(f',h')$ in $\Set$.

  Using an equalizer $E\to I\rightrightarrows J$ in $\Set$, we have
  an equalizer $G\times E\to G\times I\rightrightarrows G\times J$ in $\Setaf{G}$.
  We shall show that this (co)equalizer is preserved by $H_*$, i.e.\ the
  following diagram is a coequalizer in $\Seta{G^\op}$.
  \[ G^E \ot G^I \leftleftarrows G^J \]
  The underlying sets form a coequalizer diagram in $\Set$
  because the functor $\abs{G}^{(\place)}\colon\Set^\op\to\Set$ preserves
  reflexive coequalizers for $\abs{G}\ge 2$.
  Hence, the diagram is also a coequalizer in $\Set{G^\op}$.
\end{prf}
\begin{lemma}
  The functor $H_*\colon(\Setaf{G})^\op\to \Setaf{G^\op}$ reflects isomorphisms.
\end{lemma}
\begin{prf}{} 
  Let $f\colon X\to Y$ be a morphism $f$ in $\Setaf{G}$ such that
  $H_*(f)\colon H_*Y\to H_*X$ is an isomorphism.
  There are three cases $1\to 1$, $G\times I\to 1$, and $G\times I\to G\times J$
  for the $G$-map $f$, and the first two cases are trivial.
  For the last case, we may assume $f=\id_G\times k$ by Lem.~\ref{lemma:repnFree}.
  The right $G$-bijection $H_*(G\times J)\to H_*(G\times I)$ is
  $G^k\colon G^J\to G^I$.
  By $\abs{G}\ge 2$, the map $k$ is a bijection, which shows that
  the $G$-map $f=\id_G\times k$ is an isomorphism.
\end{prf}

\begin{prf}{ of  Lem.~\ref{lemma:posetAdjIncl}}
  For a presheaf $\lft l\in\Do P$,
  we shall show $H^*\lft l\in\Upper P$.
  The set $(H^*\lft l)(x) = (\Do P)(\lft l,\nabla x)$
  has at most one element
  for any $x\in P$,
  since $\nabla x$ is a subobject of a terminal object $1\in\Do P$.
  In particular, the postsheaf $H^*\lft l\in\Up P$
  is an upper set of $P$.

  A presheaf $\lft l$ can be written as a canonical colimit:
  $\lft l = \supp_i \nabla x_i$.
  Then, $H^*\lft l = \supp_i \nabla x_i$.
  The colimit in $\Up P$ is a limit in $\Set^P$,
  moreover it is just a product in $\Set^P$ because
  the objects $\nabla x_i$ are subobjects of $1$ in $\Set^P$.
  Let $L\subseteq P$ be the lower set
  defined by the following coproduct in $\Lower P$.
  \[ L=\bigcup_i \nabla x_i
  = \bigcup_i \{y\mid y\le x_i\} = \{y\mid \lft l(y) \ne 0\} \]
  It holds that
  $H^*L=\prod_i \nabla x_i = \supp_i \nabla x_i = H^*\lft l$.
\end{prf}

\begin{prf}{ of  Cor.~\ref{corollary:retrImPoset}}
  A retract of an upper set is also an upper set.
  \begin{align*}
    \Retr_{\Up P}(\Img H^*)
    &= \Retr_{\Up P}(\Img (H^*|_{\Lower P})) &&\text{by the above lemma}
    \\&=\Img (H^*|_{\Lower P}) &&\text{by }\Img (H^*|_{\Lower P})\subseteq\Upper P
  \end{align*}
  \vskip-\baselineskip
\end{prf}

\begin{prf}{ of  Lem.~\ref{lemma:reflEqZpZ}}
  Let $r$ be a common retraction, and
  $I=\Img f \cup \Img g \cong \od{I_1}{I_p}$.
  By the existence of retraction, we have
  $f=f'_1+f'_p\colon \od{U_1}{U_p} \to \od{U'_1}{U'_p}$
  for some $\zpz$-maps $f'_1,f'_p$,
  and similar for $g$.
  Hence, there exists $r|_I=r'_1+r'_p\colon\od{I_1}{I_p}\to\od{U_1}{U_p}$.
  We may assume $r'_1+r'_p=\od{r_1}{r_p}$ by modifying
  the coercing isomorphism $U'\cong\od{U'_1}{U'_p}$ in $I$.
  Under the assumption, we obtain
  $f'_1+f'_p=\od{f_1}{f_p}$, $g'_1+g'_p=\od{g_1}{g_p}$.
  
  The reflexive equalizer in $\Seta{\zpz}$
  is also a reflexive equalizer in $\Set$,
  which induces a pullback of injections
  \begin{center}\begin{tikzar}{}
    E_1+pE_p \arrow[>->]{d}[swap]{e_1} \arrow[>->]{r}{e_1}\& U_1+pU_p \arrow[>->]{d}{f_1}
  \\
  U_1+pU_p \arrow[>->]{r}[swap]{g_1}\& U'_1+pU'_p
  \rlap{\enspace.}
  \end{tikzar}\end{center}
  Changing the base by maps
  $U'_1\to U'_1+pU'_p$,
  $U'_p\to U'_1+pU'_p$
  concludes the proof.
\end{prf}

\begin{prf}{ of  Prop.~\ref{proposition:monadicZpZ}}
  It is easy to check the full subcategories $\CCC,\DDD$ contain
  all retracts of images.
  Then, by Lem.~\ref{lemma:retrAdj}.1,
  we have only to show that the restrictions
  $\DDD\to\CCC$, $\CCC^\op\to\DDD^\op$ are monadic.

  By the above lemma, an equalizer in $\Seta{\zpz}$
  of a reflexive pair in $\DDD^\op$ can be taken as
  \begin{center}
    \begin{tikzar}{}
  \od{E_1}{E_p}
  \arrow[->]{r}{\od{e_1}{e_p}}\&
  \od{U_1}{U_p}
  \arrow[->,shift left=3pt]{r}{\od{f_1}{f_p}}
  \arrow[->,shift right=3pt]{r}[swap]{\od{g_1}{g_p}}
  \&
  \od{U'_1}{U'_p}
    \end{tikzar}
  \quad for some equalizers $
  \begin{aligned}
  & E_1 \tto{e_1} U_1\overset{f_1}{\underset{g_1}\rightrightarrows} U'_1
  ,\\
  & E_p \tto{e_p} U_p\overset{f_p}{\underset{g_p}\rightrightarrows} U'_p
  .
\end{aligned}
  $
  \end{center}
  It is easy to show that $\od{E_1}{E_p}\in\DDD$.
  For example, if $\Phi_1=1$ then $U_1=U'_1=1$, which implies $E_1=1$.

  By Lem.~\ref{lemma:weakReflEq}.2, the diagrams
  \[
  \Phi_1^{E_1} \ot \Phi_1^{U_1} \leftleftarrows \Phi_1^{U'_1}
  \enspace,\qquad
  (\Phi_1+p\Phi_p)^{E_p} \ot (\Phi_1+p\Phi_p)^{U_p} \leftleftarrows (\Phi_1+p\Phi_p)^{U'_p}
  \]
  are split coequalizers.
  Hence, their pointwise product
  $\Phi_*(\od{}{E_p}) \ot \Phi_*(\od{}{U_p}) \leftleftarrows \Phi_*(\od{}{U'_p})$ is a (split) coequalizer.

  Let $f\colon U\to U'$ be a $\zpz$-map.
  We may assume $f$ is of the form
  \[
  \od{U_1}{(U_p^\mathrm{L}+U_p^\mathrm{R})}
  \to \od{U'}{U'_p} \]
  induced by maps
  \[
  f_1\colon U_1\to U'_1
  \enspace,\quad
  g\colon U_p^\mathrm{L}\to U'_1
  \enspace,\quad
  f_p\colon U_p^\mathrm{R}\to U'_p
  \]
  up to isomorphisms. 
  We have to prove $f_1,f_p$ are bijections and $U_p^\mathrm{L}=0$,
  if \[ \Phi_*(f)\colon
  \Phi_1^{U'_1} (\Phi_1+p\Phi_p)^{U'_p}\to
  \Phi_1^{U_1} (\Phi_1+p\Phi_p)^{U_p^\mathrm{L}} (\Phi_1+p\Phi_p)^{U_p^\mathrm{R}} \]
  is a bijection.
  It is easy to show that the map $\Phi_*(f)$ factors through
  $\Phi_1^{U_1} \Phi_1^{U_p^\mathrm{L}} (\Phi_1+p\Phi_p)^{U_p^\mathrm{R}}$.
  Therefore, we have
  \[
  \Phi_1^{U_1} (\Phi_1+p\Phi_p)^{U_p^\mathrm{R}} = 0
  \quad\text{or}\quad
  \Phi_p = 0
  \quad\text{or}\quad
  U_p^\mathrm{L} = 0
  \enspace.
  \]
  The rest is straightforward.
  For $\Phi_p=0$ cases,
  we remark that $U_p^\mathrm{L}+U_p^\mathrm{R}=0$
  and then the claim reduces to the monadicity of
  a restriction of $\Phi_1^{(\place)}\colon\Set^\op\to\Set$.

  By the monadicity theorem, the restriction $\Phi_*\colon \DDD\to\CCC$
  is monadic.
  The other monadicity is easy.
\end{prf}

\section{Appendix: General propositions}
\begin{proposition}[Precise monadicity theorem]
  Let $U\colon\DDD\to\CCC$ be a functor
  that has a left adjoint $F\colon\CCC\to\DDD$,
  and $T=U\circ F$ be the induced monad.
  \begin{center}
    \begin{tikzar}{}
      \CCC
      \arrow[->,loop above]{}{T}
      \arrow[->,bend right=15]{d}[swap]{F} \arrow[<-,bend left=15]{d}{U} \arrow[no line]{d}{\scriptstyle\dashv}
      \arrow[->,bend right=15]{r}[swap]{} \arrow[<-,bend left=15]{r}{} \arrow[no line]{r}{\scriptstyle\top}
      \&
      \CCC^T
      \\
      \DDD
      \arrow[->,shift right=5pt]{ur}[swap]{K}
      \arrow[<-,dashed,shift left=0pt]{ur}{L}
    \end{tikzar}
$\begin{aligned}
  KD &= (UD \oot{U\varepsilon D} UFUD)
  \enspace,\\
  L(C\oot{h} UFC) &\ot FC \overset{Fh}{\underset{\varepsilon FC}\leftleftarrows} FUFC
  \text{\quad is a coequalizer.}
\end{aligned}$
  \end{center}
  \begin{enumerate}[(a)]
    \item 
      The comparison functor $K\colon\DDD\to\CCC^{T}$
      has a left adjoint $L\colon \CCC^{T}\to\DDD$ if
      the category $\DDD$ has reflexive $U$-split coequalizers.
    \item 
      The functor $L$ is full and faithful
      if $\DDD$ has and $U$ preserves reflexive $U$-split coequalizers.
    \item       The comparison functor $K$ is full and faithful if
      $U$ reflects isomorphisms \cite[Sec.~3.3]{BarrM:TTT}
 .
  \end{enumerate}
  In particular, the right adjoint functor $U$ is monadic if
  $U$ creates reflexive $U$-split coequalizers.

  Conversely, for a monad $T$,
  the forgetful functor $U^T\colon\CCC^T\to\CCC$ creates
  $U^T$-split coequalizers.
\end{proposition}

\begin{lemma}\label{lemma:weakReflEq}
  Let $X,Y$ be sets, and $f,g\colon X\rightrightarrows Y$ be a pair of maps.
  Let $E \tto{e} X$ be an equalizer of the pair $(f,g)$.
  \begin{enumerate}
    \item
      If the pair $(f,g)$ is reflexive,
      the diagram below is a pullback and
      the maps $f,g$ are injections.
    \item
      Let $R$ be a nonempty set.
      If the diagram below is a pullback and the map $f$ is an injection,
      $R^E \oot{R^e} R^X \overset{R^f}{\underset{R^g}\leftleftarrows} R^Y$
      is a split coequalizer.
      In particular, the functor $R^{(\place)}\colon\Set^\op\to\Set$ preserves
      such coequalizers.
  \end{enumerate}
  \begin{center}\begin{tikzar}{}
    E \arrow[->]{d}[swap]{e} \arrow[->]{r}{e} \& X \arrow[->]{d}{f} \\
    X \arrow[->]{r}[swap]{g} \& Y
  \end{tikzar}\end{center}
\end{lemma}
\begin{prf}{ of  Lem.~\ref{lemma:weakReflEq}}
  The item~1.\ is easy.

  Let $r\in R$ be an element.
  By injectivity, define maps
  $R^E \tto{e_r} R^X \tto{f_r} R^Y$ by
  \[ f_r(h)(y) =
  \{\; h(x)\text{\quad if }y = f(x) ;\qquad r\text{\quad otherwise} \;\}
  \quad\text{for }h\colon X\to R \enspace,\]
  and similarly for $e_r$.
  The maps gives a splitting of the diagram
      $R^E \oot{R^e} R^X \overset{R^f}{\underset{R^g}\leftleftarrows} R^Y$,
  which concludes the proof of the item.~2.
\end{prf}

\end{document}